\documentclass[leqno,a4paper,12pt]{article}

\usepackage[utf8]{inputenc}
\usepackage[LGR,T1]{fontenc}
\usepackage[greek,francais]{babel}
\usepackage{amsmath,amssymb,amscd}
\usepackage{pst-tree,pst-node}
\usepackage{ifthen}
\usepackage{enumerate}

\title{La théorie arithmétique des grandeurs algébriques de Kronecker (1882)
\footnote{Mots clefs~: Kronecker, algebraic geometry, algebraic number theory, elimination theory, 01A55, 13P15, 14-03, 11-03. Ce texte développe et corrige un chapitre de ma thèse~; il a fait l'objet d'un exposé à la << Journée en l'honneur de Christian Houzel >> à l'IHP (Paris) le 23 novembre 2007.}}
\author{Erwan Pench\`evre}
\date{}

\begin{document}

\newcommand{\pagesname}{p.}
\renewcommand{\pagename}{p.}
\newcommand{\Inname}{\textit{in}}
\newcommand{\volumename}{vol.}
\newcommand{\ofname}{,}
\newcommand{\editorname}{ed.}
\newcommand{\technicalreportname}{Th\`ese}
\renewcommand{\up}[1]{\textsuperscript{#1}}
\newcommand{\gche}{ << }
\newcommand{\drte}{ >> }
\newcommand{\und}[2]{\underset{#1}{#2}}
\renewcommand{\bf}[1]{\mathbf{#1}}
\newcommand{\esprit}{`}
\newcommand{\itnor}[1]{\textit{#1}}
\newcommand{\ul}[1]{\"{#1}}
\newcommand{\biblio}[1]{[#1]}
\newcommand{\bino}[2]{{{#1}\choose{#2}}}
\newcommand{\under}[1]{\underline{#1}}
\newcommand{\im}{\text{Im}}
\newcommand{\coker}{\text{coker}}
\newcommand{\ord}{\text{ord}}
\newcommand{\card}{\text{Card}}
\newcommand{\cf}{\textit{cf}. }
\newcommand{\Cf}{\textit{Cf}. }
\newcommand{\ibid}{\textit{ibid}. }
\newcommand{\limproj}{\underleftarrow{\text{lim}}}
\newcommand{\eqcart}{\reflectbox{$\varpropto$}}
\newcommand{\diff}{\text{d}}
\newcommand{\codim}{\text{codim}}

\newtheorem{theo}{Théorème}
\newtheorem{lemme}{Lemme}
\newtheorem{corol}{Corollaire}
\newcommand{\lc}[1]{\text{lc}#1}
\renewcommand{\sc}[1]{\text{sc}#1}
\newboolean{note}
\setboolean{note}{false}  
\newcommand{\note}[1]{\ifthenelse{\boolean{note}}{\textbf{[#1]}}{}}

\newcommand{\rat}{\mathfrak{R}}
\newcommand{\gat}{\mathfrak{G}}
\newcommand{\spe}{\mathfrak{S}}
\newcommand{\fon}{\mathfrak{f}}
\newcommand{\Fm}{\text{Fm}}
\newcommand{\Nm}{\text{Nm}}

\newcommand{\et}{\text{\&}}
\newcommand{\etc}{\text{\&c. }}

\maketitle

\label{Kronecker}
La préface des \textit{Grundz\ul{u}ge einer arithmetischen Theorie der algebraischen Gr\ul{o}ssen} annonce en quelques phrases le projet de Kronecker, qui, selon Dieudonné (\cite{dieudonne1974}, p.~27), <<~was the first to dream of one vast algebraico-geometric construction comprising the theory of numbers and algebraic geometry~>>. 
La première partie des \textit{Grundz\ul{u}ge} est une théorie des extensions de corps et des entiers algébriques. Comme l'indique Kronecker%
\footnote{\Cf p. 2 où il parle du problème de <<~la représentation de toutes les grandeurs entières algébriques d'un genre... par des fonctions linéaires de plusieurs grandeurs algébriques... afin de rendre réelles les grandeurs fractionnaires \textit{idéales}~>>} et comme l'a démontré H. M. Edwards, une bonne définition du concept même d'entier est à cet égard cruciale. Dans la deuxième partie, il s'agit alors de faire une théorie de la divisibilité dans les anneaux d'entiers. Mais, comme le remarque J. Boniface, à la différence de Dedekind, Kronecker, lors du passage des entiers naturels aux entiers algébriques, ne s'intéresse pas tant à préserver les lois (théorème d'existence et d'unicité de la décomposition en facteurs premiers, qui joue un rôle majeur chez Dedekind, mais est secondaire, comme l'a montré Edwards, chez Kronecker) qu'à préserver les objets, ici le plus grand commun diviseur de plusieurs entiers, qu'il définit pour un anneau d'entiers quelconque. Enfin, la théorie des corps construite par Kronecker concerne aussi bien les extensions transcendantes, qu'il manipule géométriquement au moyen d'un concept de <<~variété~>> (\textit{Mannigfaltigkeit}), et il laisse entendre à la fin de sa préface que le concept de <<~Stufe~>> (ce qui, traduit dans un langage algébrico-géométrique plus moderne, désigne la codimension d'une sous-variété de l'espace affine) aura un rôle important dans l'édifice%
\footnote{Selon Kronecker, à l'issue du dernier paragraphe des \textit{Grundz\ul{u}ge}, il apparaîtra que <<~le concept de \textit{Stufe} convient pour remplacer le concept d'irrationalité algébrique~>>. Nous y reviendrons \textit{infra} section \ref{Stufe}.}. %
Nous verrons que ce concept est au c{\oe}ur de la théorie de l'élimination élaborée par Kronecker.

\section{Théorie des corps}
Kronecker choisit délibéremment (il s'en justifie dans le \S~1) d'éviter le mot <<~corps~>> utilisé par Dedekind, et il préfère parler de <<~domaine de rationalité~>> (\textit{Rationalit\ul{a}ts-Bereich}) pour désigner un corps $\mathbb{Q}(\rat',\rat'',\rat''',...)$ qu'il note lui-même de la manière suivante~:
$$(\rat',\rat'',\rat''',...)$$
où $\rat'$, $\rat''$, $\rat'''$,... sont des <<~grandeurs~>> données. Kronecker insiste d'emblée sur le fait que le concept de <<~grandeur~>> recouvre ici les notions usuelles de nombres rationnels et nombres irrationnels algébriques (tel $\sqrt{2}$), aussi bien que celle de <<~grandeur indéterminée~>> ou de <<~fonction rationnelle de plusieurs grandeurs indéterminées~>>. Pour motiver sa démarche, il indique dès le \S~1 que la notion d'irréductibilité d'un polynôme (à coefficients dans un domaine de rationalité) n'a de sens que relativement au domaine de rationalité choisi, et il remarque que Galois et Abel étaient déjà conscients du besoin de préciser le sens du mot <<~rationnel~>>. Il renvoie au manuscrit d'Abel \textit{Sur la résolution algébrique des équations} publié après sa mort en 1839. Abel y définit en effet les différentes modalités du concept de <<~grandeur~>> presque comme Kronecker, et il écrit~:
\begin{quote}
D'après la nature des quantités connues nous ferons les distinctions suivantes~:
\begin{enumerate}[1.]
\item Une quantité qui peut s'exprimer algébriquement par l'unité s'appelle un nombre algébrique~; si elle peut s'exprimer rationnellement par l'unité, elle s'appelle un nombre rationnel, et si elle peut être formée de l'unité par addition, soustraction et multiplication, elle s'appelle un nombre entier.
\item Si les quantités connues contiennent une ou plusieurs quantités variables, la quantité $y$ est dite fonction algébrique, rationnelle ou entière de ces quantités selon la nature des opérations nécessaires pour la former. Dans ce cas on regarde comme quantité connue toute quantité constante.%
\footnote{\cf \cite{abel1881} t.~2 p.~224. Ailleurs, sur le caractère relatif de la notion d'irréductibilité, Abel s'exprimait ainsi~: <<~Une équation $\phi x=0$, dont les coefficiens sont des fonctions rationnelles d'un certain nombre de quantités connues $a$, $b$, $c$,... s'appelle \textit{irréductible}, lorsqu'il est impossible d'exprimer aucune de ses racines par une équation moins élevée, dont les coefficiens sont également des fonctions rationnelles de $a$, $b$, $c$,...~>> (\cite{abel1881}, t.~I, p.~479).}
\end{enumerate}
\end{quote}
\par
Le lexique introduit par Kronecker au \S~2 de son mémoire permet de décrire le contexte relatif, que l'on exprime aujourd'hui en précisant de quel corps de base un corps donné est l'extension. Kronecker nomme <<~domaine de base~>> (\textit{Stammbereich}) un corps de base $\mathbb{Q}(\rat',\rat'',\rat''',...)$, et <<~domaine de genre $\gat'$~>> (\textit{der Bereich der Gattung $\gat'$}) une extension engendrée par la grandeur $\gat'$, c'est-à-dire de la forme $\mathbb{Q}(\gat',\rat',\rat'',\rat''',...)$, ce qu'il note lui-même~:
$$(\gat',\rat',\rat'',\rat''',...)$$
Un genre $\gat'$ est dit <<~d'ordre $n$~>> quand $\gat'$ est de degré $n$ sur le corps de base~; il est dit <<~genre galoisien~>> quand l'extension est galoisienne (c'est-à-dire que tous les conjugués de $\gat'$ appartiennent à $\mathbb{Q}(\gat',\rat',\rat'',\rat''',...)$). Dans la terminologie de Kronecker, l'expression <<~genre $\gat'$~>> ne s'applique en fait qu'à l'ensemble des éléments primitifs (c'est-à-dire les éléments de même degré que $\gat'$ sur le corps de base) de l'extension considérée, tandis que l'expression <<~domaine de genre $\gat'$~>> désigne toute l'extension. Kronecker démontre que le degré d'une sous-extension divise le degré de l'extension.
\par
Au \S~3, Kronecker remarque que l'on peut construire toute extension de $\mathbb{Q}$ par l'adjonction <<~d'un certain nombre de grandeurs variables ou indépendantes, et de fonctions algébriques de celles-ci~>>. Il démontre ensuite le théorème de l'élément primitif, qui permet d'écrire toute extension de type fini de $\mathbb{Q}$ sous la forme $\mathbb{Q}(\rat',\rat'',\rat''',...)$ où tous les $\rat$ sauf un sont des indéterminées. Il mentionne aussi l'existence d'extensions qui ne sont pas de type fini, telle la clôture algébrique de $\mathbb{Q}$, mais il choisit de restreindre son propos aux extensions de type fini%
\footnote{En fait, Kronecker n'impose pas de limite au nombre d'interminées $\rat',\rat'',\rat'',...$. Il pourrait y en avoir une infinité, et il serait donc plus exact de parler, de manière générale, d'<<~extensions finies d'extensions purement transcendantes de $\mathbb{Q}$~>>.}.
\par
Par anticipation, remarquons déjà que tout le développement des \textit{Grundzüge} se rapportera tantôt à l'une, tantôt à l'autre des trois structures suivantes~:
\begin{enumerate}[--]
\item un anneau d'entiers $O_K\subset K$ où $K$ est une extension algébrique de $\mathbb{Q}$ (c'est l'objet de la théorie des nombres algébriques)
\item un anneau de degré de transcendance $(n-1)$ sur $\mathbb{Q}$, engendré par $(n-1)$ indéterminées~: $\mathbb{Q}[x_1,...,x_{n-1}]\subset\mathbb{Q}(x_1,...,x_{n-1})$ (c'est l'objet de la théorie des variétés affines sur un corps, et de la théorie de l'élimination)
\item un anneau d'entiers $O_K\subset K$ où $K$ est une extension algébrique $k(x_n)$ d'un corps $k=\mathbb{Q}(x_1,...,x_{n-1})$ de degré de transcendance $(n-1)$ sur $\mathbb{Q}$
\end{enumerate}
Mais quand on parle de <<~structures algébriques~>> dans l'{\oe}uvre de Kronecker, il faut bien-sûr y mettre un bémol puisqu'il avait un point de vue finitiste et refusait l'usage d'ensembles infinis.
\par
Dans cet ordre d'idée, au \S~4 apparaît une modalité importante de la pensée de Kronecker~: l'exigence d'effectivité. La théorie des corps qu'il vient d'exposer permet certes une définition rigoureuse de la notion d'irréductibilité, relative à un corps de base, mais cette définition reste incomplète, selon Kronecker, tant qu'elle n'est pas accompagnée d'un algorithme pour tester l'irréductibilité d'un polynôme donné. Plus précisément, il exige un algorithme pour tester l'irréductibilité d'un polynôme $F(x,x',x'',...)$ à coefficients dans une extension de type fini de $\mathbb{Q}$. Il remarque d'emblée que d'habitude on avait recours à la théorie de l'élimination pour résoudre ce problème. Et que l'on pense en effet aux travaux des algébristes du XVII\up{ème} pour la recherche de facteurs rationnels (Jan Hudde) par la méthode des coefficients indéterminés, puis aux travaux de Lagrange que nous avons déjà commentés. Les raisonnements sur lesquels Bézout (et plus tard Cayley) fonde sa théorie de l'élimination, et que les progrès de l'algèbre linéaire jusqu'à Kronecker ne peuvent éclaircir d'avantage, ont un défaut que la théorie de l'élimination conçue par Poisson (puis Hesse et Schläfli) ne peut relever, puisque celle-ci se place d'emblée dans le cadre restreint du résultant d'un système de $n$ équations ayant un nombre fini de solutions. Une théorie plus générale était encore à naître. Revenons au \S~4 des \textit{Grundzüge}. Kronecker juge que~:
\begin{quote}
...lors du développement naturel et complet de la théorie de l'élimination on utilise la décomposition d'une fonction entière en ses facteurs. C'est pourquoi on aura recours ici à une nouvelle méthode [de décomposition] qui ne nécessite que des moyens simples, déjà applicable ici.
\end{quote}
Ce qui semblait d'abord n'être qu'une exigence gratuite d'effectivité revêt donc une autre dimension. Kronecker entend fonder la théorie de l'élimination. Certes la théorie kroneckerienne de l'élimination est une théorie effective~: il ne s'agit ni de donner une définition \textit{a priori} du résultant, ni de donner seulement une borne au nombre de solutions d'un système d'équations polynomiales, mais vraiment de résoudre un tel système absolument quelconque. La méthode d'élimination proposée par Kronecker consiste en fait, interprétée en termes modernes, en la décomposition d'une variété algébrique affine (ensemble des solutions d'un système d'équations polynomiales) en ses composantes irréductibles. Et cette méthode requiert en effet un algorithme pour décomposer un polynôme en ses facteurs irréductibles. Décrivons rapidement celui que donne Kronecker. En faisant le produit des conjugués du polynôme $F(x,x',x'',...)$ et en posant 
$$x'=c_1x^g,\quad x''=c_2x^{g^2},...\quad x^{(n)}=c_nx^{g^n}$$
pour $g$ assez grand, il se ramène au cas d'un polynôme $F(x)$ à une indéterminée, à coefficients dans $\mathbb{Z}$. Dans ce cas, il représente le facteur $f(x)$ cherché par une formule d'interpolation
$$f(x)=f(r_0)g_0(x)+f(r_1)g_1(x)+...+f(r_n)g_n(x)$$
où les $r_i$ sont des entiers distincts quelconques, et les $g_n(x)$ certains polynômes. Si $f(x)$ est un facteur de $F(x)$, alors pour tout $i$, $f(r_i)$ divise $F(r_i)$, et il n'y  a donc qu'un nombre fini de facteurs $f(x)$ possibles.
A ce stade plane un mystère~: quelle place peut bien avoir la théorie de l'élimination au sein d'une <<~théorie arithmétique des grandeurs algébriques~>>~? Une première réponse nous sera donnée par le rôle qu'occupe le discriminant dans la théorie arithmétique de Kronecker.

\section{Grandeur entière et discriminant}
Au \S\S~5--7, Kronecker définit la <<~grandeur entière~>>. Il utilise à cet effet une nouvelle notation 
$$[\rat',\rat'',\rat''',...]$$
pour désigner l'anneau que l'on noterait aujourd'hui $\mathbb{Z}[\rat',\rat'',\rat''',...]$. Kronecker classe les grandeurs entières en <<~espèces~>> (\textit{Art}) au sein d'un <<~genre~>>. Si $\rat'$, $\rat''$, $\rat'''$,... sont des indéterminées, les entiers de $\mathbb{Q}(\rat',\rat'',\rat''',...)$ sont simplement les éléments de $\mathbb{Z}[\rat',\rat'',\rat''',...]$%
\footnote{Une définition de la notion d'entier algébrique commence nécessairement par la donnée, assez arbitraire, d'un anneau d'entiers en quelque sorte <<~absolu~>>, intégralement clos, dont chaque autre anneau d'entiers sera une certaine clôture intégrale dans une extension de son corps de fractions. En caractéristique zéro, le choix naturel est celui de $\mathbb{Z}$. Mais parfois, Kronecker énonce des résultats <<~en faisant abstraction des coefficients entiers~>> (\textit{sofern von den Zahlcoefficienten abgesehen wird}, \cf par exemple début du \S~7). Cela revient alors à définir l'anneau d'entiers absolu comme étant $\mathbb{Q}[\rat',\rat'',\rat''',...]$, où $\rat'$, $\rat''$, $\rat'''$,... sont des indéterminées (remarquez que dans ces conditions, aussi bien $\mathbb{Z}[\rat',\rat'',\rat''',...]$ que $\mathbb{Q}[\rat',\rat'',\rat''',...]$ sont intégralement clos dans leur corps de fraction commun $\mathbb{Q}(\rat',\rat'',\rat''',...)$).}. %
Dans une extension de corps, une grandeur est dite entière si elle est solution d'une équation à coefficients entiers dans le corps de base et coefficient dominant égal à 1. Si $\spe'$, $\spe''$, $\spe'''$,... sont des entiers algébriques dans un <<~domaine de genre~>> $K=\mathbb{Q}(\gat,\rat',\rat'',\rat''',...)$ sur un corps de base $k=\mathbb{Q}(\rat',\rat'',\rat''',...)$, Kronecker nomme l'anneau $\mathbb{Z}[\rat',\rat'',\rat''',...,\spe',\spe'',\spe''',...]$ un <<~domaine d'espèce~>>. La réunion de tous les domaines d'espèces d'un même domaine de genre constitue la clôture intégrale de l'anneau d'entiers du corps de base dans son extension. L'objectif des \S\S~6 et 7 est de démontrer que cet ensemble est bien lui-même un <<~domaine d'espèce~>>, que Kronecker qualifiera d'<<~espèce principale~>> (\textit{Haupt-Art}). Kronecker démontre un peu plus~: notons $O_k$ l'anneau d'entiers du corps de base $k$ et $O_K$ sa clôture intégrale dans $K$, alors $O_K$ est un $O_k$-module de type fini. Au \S~6, il montre que $O_K$ est un $O_k$-module de type fini, en se ramenant au cas où $k$ est une extension purement transcendante de $\mathbb{Q}$ (au \S~7, il remarque que dans certains cas, un nombre $n$ de générateurs suffisent, où $n$ est le degré de $K$ sur $k$)%
\footnote{A partir du \S~6, et au moins jusqu'au \S~18, Kronecker utilise souvent l'hypothèse que $k$ est une extension purement transcendante de $\mathbb{Q}$ (\textit{natürliche Rationalitäts-Bereich}). C'est le fait que $k$ ait un anneau d'entiers factoriel qui importe dans les démonstrations. On n'obtient donc une théorie relative de la divisibilité dans les anneaux d'entiers que si l'anneau d'entiers du corps de base a lui-même <<~naturellement~>> une telle théorie. De toute façon, comme on l'a déjà remarqué, tous les corps considérés par Kronecker sont des extensions finies d'extensions purement transcendantes de $\mathbb{Q}$.}. %
A cet effet, il introduit le concept de \textit{discriminant de $n$ grandeurs} $x^{(1)},x^{(2)},...,x^{(n)}$ d'un genre d'ordre $n$~: c'est le carré du déterminant $\left\vert x_i^{(h)}\right\vert$ ($h,i=1,2,...,n$) où pour tout $h$, $x_1^{(h)}$, $x_2^{(h)}$,..., $x_n^{(h)}$ sont les conjugués%
\footnote{Lorsque les $n$ grandeurs sont des entiers algébriques, le fait que les conjugués d'un entier soient eux-mêmes des entiers implique que le discriminant est un entier algébrique.} %
de $x^{(h)}$.
Ce concept était déjà présent, sans être nommé, dans son article de 1862 (publié en 1881) \textit{\ul{U}ber die Discriminante algebraischer Functionen einer Variabeln}. Comme son titre l'indique, ce texte se restreignait au cas des grandeurs algébriques sur un corps de base de la forme $k=\mathbb{Q}(v)$ (où $v$ est une indéterminée). Kronecker y étudiait le \textit{discriminant d'une grandeur $x$} solution d'une équation irréductible de degré $n$ à coefficients dans ce corps, discriminant défini comme étant le carré de $\left\vert x_i^{h}\right\vert$ où les $x_i^{h}$ (pour $1\leq i\leq n$ et $0\leq h\leq n-1$) sont les conjugués des $n$ premières puissances de $x$. Le \textit{discriminant d'une grandeur $x$} est donc le \textit{discriminant des $n$ grandeurs $1,x,...,x^{n-1}$}, et il coïncide d'ailleurs avec le \textit{discriminant de l'équation}, défini comme produit des différences des racines. Soit $x$ une grandeur entière algébrique sur $k$, et $K=k(x)$. Dans ce cas particulier, $n$ générateurs suffisent à engendrer le $O_k$-module $O_K$. Kronecker montre que, pour tout autre grandeur algébrique entière $y$ sur $k$ telle que $K=k(y)$, le discriminant $\Delta$ de ces $n$ générateurs divise le discriminant de $y$. Il l'appelle alors <<~facteur essentiel~>> du discriminant de $y$ (\textit{wesentliche Theiler}). Il montre même que $\Delta$ est le plus grand commun diviseur de tous ces discriminants%
\footnote{\Cf p.~221.}, ce qui lui donne un sens intrinsèque (indépendant du choix des générateurs). Il montre de plus que le <<~facteur non essentiel~>> du discriminant de $y$ est toujours un carré. Dans la préface de 1881 à cet article, Kronecker écrivait~:
\begin{quote}
Je pourrais décrire comme étant l'un des principaux fruits de ce travail le fait que le traitement du cas pour ainsi dire général ou régulier des équations algébriques entre deux variables, auquel Riemann s'est restreint en remarquant que <<~l'on peut aisément étendre les résultats aux autres cas, regardés comme cas limites~>>, suffit aussi parce que les autres équations peuvent être transformées en de telles.%
\footnote{\Cf p.~199.}
\end{quote}
Riemann, dans sa \textit{Théorie des fonctions abéliennes}, s'était en effet restreint au cas des fonctions dont les points de ramification <<~coïncident seulement par paires, et en se détruisant~>> (\cf \S~VI). Kronecker démontre qu'il existe un $\xi\in K$ dans le discriminant duquel le facteur essentiel et le facteur non essentiel n'ont aucune racine commune et les racines du facteur non essentiel sont toutes des racines de multiplicité exactement égale à 2. Géométriquement, les racines du facteur essentiel sont le lieu de ramification de $\xi$, vu comme fonction de la variable $v$, et les racines du facteur non essentiel correspondent à ce que Riemann appelait des <<~points de ramification qui coïncident~>>, c'est-à-dire, dans un langage plus moderne, à des points singuliers de la courbe d'équation $F(\xi,v)=0$ (remarquons que la <<~surface de Riemann~>> de $\xi$ n'est autre qu'un modèle non singulier de cette courbe, au sein de sa classe d'équivalence birationnelle). Kronecker est donc parvenu à transformer l'équation en une équation dont les points singuliers sont tous des points doubles d'abscisses distinctes et situés hors de son lieu de ramification. Noether (\cf \cite{brill1892}, p.~370--375) y a vu la première méthode de résolution des singularités des courbes algébriques planes, préfigurant les résultats de peu postérieurs (début des années 1870) de Noether, Hamburger et Weierstrass (\cf \cite{penchevre2006c} section 13.2). 
\par
Mais revenons aux \textit{Grundz\ul{u}ge}. Le discriminant (et le lieu de ramification) était un invariant important en théorie des nombres (Dirichlet avait par exemple donné une formule dépendant du discriminant pour le nombre de classes d'un corps quadratique) ainsi qu'en théorie des fonctions d'une variable (par exemple pour déterminer le genre d'une surface de Riemann). Kronecker cherche à présent (\S~8) à généraliser cette notion au cas de plusieurs variables. Mais, dans le cas général, le système de générateurs du $O_k$-module $O_K$ a plus de $n$ éléments, et Kronecker considère le <<~système fondamental de discriminants~>>, ensemble des discriminants des générateurs pris $n$ à $n$, ainsi que (dans le langage d'aujourd'hui) l'idéal engendré par tous ces discriminants\label{KroneckerIdealDeterminantiel}%
\footnote{Dans la terminologie moderne, cet idéal est un <<~idéal déterminantiel~>>, \cf \cite{giusti1981}.}. %
Géométriquement, c'est-à-dire quand les indéterminées du corps de base sont conçues comme des variables, il explique que la variété des zéros de cet idéal est un <<~invariant du genre~>>. Dans son langage, il représente en fait cet idéal par une forme linéaire dont les coefficients sont les discriminants des générateurs pris $n$ à $n$ (on verra plus loin le rapport avec le concept d'idéal). Dans le cas où l'anneau d'entiers du corps de base est factoriel, il propose de définir le \textit{discriminant du genre} comme étant le plus grand commun diviseur $\Delta$ des discriminants des générateurs pris $n$ à $n$ (sa variété des zéros est, nous dit Kronecker, une partie de la variété des zéros de l'idéal ci-dessus). Kronecker montre qu'alors, si $x'$, $x''$,..., $x^{(n+m)}$ sont les générateurs de $O_K$, et $u_1x'+u_2x''+...+u_{n+m}x^{(n+m)}$ un élément générique de $O_K$, le facteur de son discriminant%
\footnote{On remarquera de plus que la non-annulation du discriminant de $u_1x'+u_2x''+...+u_{n+m}x^{(n+m)}$, vu comme polynôme en $u_1,...,u_n$, est une condition pour que cet élément soit un élément primitif (puisqu'alors tous ses conjugués sont distincts). Kronecker utilise ainsi le discriminant dans un cas particulier, fin \S~12.} %
indépendant des $u$ divise une puissance de $\Delta$. 
\par
Jusqu'ici, Kronecker a donc posé les fondements d'une théorie des extensions de corps et des anneaux d'entiers, dans laquelle le polynôme minimal d'une grandeur (qui détermine ses conjugués, l'ordre du genre auquel elle appartient, son discriminant, et est unitaire à coefficients entiers si et seulement si la grandeur est entière%
\footnote{On est ici en caractéristique 0, donc si $\alpha$ est un entier algébrique, son polynôme minimal est $\prod(X-\alpha_i)$, où les $\alpha_i$ sont les conjugués distincts de $\alpha$, entiers algébriques eux aussi. L'équation $\prod(X-\alpha_i)=0$ est alors une équation de dépendance intégrale de $\alpha$.}) joue un rôle important. Mais une théorie complète de la <<~grandeur algébrique~>> ne peut faire l'économie de l'étude des <<~systèmes absolument quelconques d'équations~>> (\textit{ganz allgemeiner Gleichungssysteme}). Il faudrait par exemple généraliser le concept de discriminant d'une grandeur algébrique (solution d'une équation à coefficients dans le corps de base) au cas d'un système de grandeurs, solution d'un système d'équations à plusieurs inconnues. Au \S~10, Kronecker remarque que les rapports des $(n+1)$ grandeurs $x^0, x',..., x^{(n)}$ solutions d'un système de $n$ équations homogènes $F_1=0, F_2=0,..., F_n=0$ sont algébriques et définissent un <<~genre~>> (en toute rigueur, il faut supposer ici qu'il n'y a qu'un nombre fini de solutions)~: cet ensemble de grandeurs engendre une extension algébrique du corps des coefficients des équations. Dans ce cas, il propose de déterminer directement le discriminant du genre, au moyen d'un jacobien (généralisant ainsi le discriminant d'une équation à une inconnue, résultant de l'équation et de sa dérivée). Il se donne une $(n+1)$-ième forme de degré quelconque $F_0$ à coefficients indéterminés et considère le résultant des $(n+1)$ équations
$$F_1=0, F_2=0,..., F_n=0, \left\vert F_{gh}\right\vert=0$$
où $\left\vert F_{gh}\right\vert$ est le jacobien de $F_0, F_1,..., F_n$. Il affirme alors, sans le démontrer, que le discriminant du genre%
\footnote{Pour donner un sens à cette affirmation, il faut préciser quel est l'anneau d'entiers du corps de base, soit
par exemple $\mathbb{Z}[\fon',\fon'',...]$ où $\fon',\fon'',...$ sont les coefficients (génériques) des équations. Pour une expression plus classique de ce <<~discriminant~>>, \cf \cite{molk1992}, tome I, volume 2, \S~66, p.~147--148.} %
est un facteur irréductible (quand les $F$ sont génériques) de ce résultant, indépendant des coefficients de $F_0$, et dont l'annulation est une condition nécessaire et suffisante pour que deux des solutions du système d'équations initial soient égales.

\section{Théorie de l'élimination}
Kronecker se propose de publier ultérieurement une étude générale sur les <<~systèmes absolument quelconques d'équations~>> et d'en indiquer seulement les principaux résultats dans les \textit{Grundz\ul{u}ge}. Il décrit alors, en un paragraphe que nous allons traduire intégralement, un véritable programme de recherches pour la théorie de l'élimination, au sein de laquelle le contexte relatif des <<~domaine de base~>> et <<~domaine de genre~>> aura désormais un rôle à jouer~:
\begin{quote}
Un système d'un nombre quelconque d'équations algébriques pour $z^0, z', z'',..., z^{(n-1)}$, dont les coefficients appartiennent au domaine de rationalité $(\rat',\rat'',\rat''',...)$, définit des relations algébriques entre les grandeurs $z$ et $\rat$, dont la connaissance et la représentation constituent le but de la théorie de l'élimination. Les $m$ fonctions $G_1, G_2,..., G_m$ qui, posées égales à zéro, constituent les équations, doivent être supposées fonctions rationnelles entières des $n$ grandeurs $z$, mais dont les coefficients sont seulement fonctions rationnelles (entières ou fractionnaires) des grandeurs $\rat$ à coefficients entiers. Quant au nombre des fonctions, que l'on note $m$, il ne faut aucunement le restreindre; il peut, comme dans le cas particulier -- dit général -- ci-dessus [cas du discriminant de $n$ équations homogènes à $(n+1)$ inconnues] étre égal à $n$, c'est-à-dire au nombre des grandeurs $z$ à déterminer, mais il peut aussi être plus grand ou plus petit que ce nombre. Si l'on conçoit les grandeurs $z$ comme des variables libres (\textit{als unbeschränkt veränderlich}), alors les équations $G=0$ constituent une certaine restriction de cette variabilité (\textit{eine gewisse Beschränkung dieser Variabilität}), dont une caractérisation plus précise peut être décrite comme étant le problème de l'élimination. Il faut alors maintenir l'indétermination des grandeurs indéterminées parmi les grandeurs $\rat', \rat'', \rat''',...$, ou leur libre variabilité si elles sont conçues comme des variables; c'est-à-dire que seule la variabilité des grandeurs $z$ est à caractériser, ce qui ne demande aucune restriction de la variabilité des grandeurs $\rat$. La distinction ainsi faite entre les variables $z$ et celles qui surviennent parmi les grandeurs $\rat$ est de la plus grande importance; elle n'implique aucune restriction de la généralité, elle sépare plutôt seulement les problèmes distincts qui peuvent être posés par l'élimination, selon leur contenu conceptuel.
\end{quote}
On dirait aujourd'hui que Kronecker se propose d'étudier un idéal définissant une variété plongée dans un espace affine $\mathbb{A}^n_{\mathbb{Q}(\rat',\rat'',\rat''',...)}$. Un des <<~principaux résultats~>> de la théorie de Kronecker réside en la possibilité de ramener le cas d'un système quelconque d'équations au cas d'une unique équation. Kronecker l'explique géométriquement à la fin du \S~10~: toute composante $Y$ de dimension $k$ de la variété $X\subset\mathbb{A}^n$ peut être projetée birationnellement%
\footnote{\label{AlternativeReecriture}Pour retrouver ces résultats dans un contexte moderne, \cf \cite{hartshorne1977}, exercice I.4.8 p.~31, et l'article récent de Marc Giusti et Joos Heintz \cite{giusti1991a}, qui étudie le coût (théorie de la complexité) de la méthode de Kronecker. Enfin \cite{giusti2001} et \cite{lecerf2001} suggèrent que la méthode de Kronecker est une alternative efficace aux algorithmes d'élimination ayant recours à des <<~techniques de réécriture~>> (c'est-à-dire aux bases de Gröbner).} %
sur une hypersurface d'une certaine sous-variété linéaire de dimension $(k+1)$.  Soient $x^{(0)},x^{(1)},...,x^{(n-1)}$ des coordonnées suffisamment générales de l'espace affine $\mathbb{A}^n$~; il suffira de projeter sur la sous-variété linéaire 
$$\lbrace x^{(k+1)}=x^{(k+2)}=...=x^{(n-1)}=0\rbrace\simeq\mathbb{A}^{k+1}.$$ 
D'ailleurs, la projection de la composante de dimension $k$ sur la sous-variété 
$$\lbrace x^{(0)}=x^{(k+1)}=x^{(k+2)}=...=x^{(n-1)}=0\rbrace\simeq\mathbb{A}^k$$
en fait un revêtement ramifié. On a en effet le diagramme commutatif suivant~:\label{NormalisationNoether}
$$\hspace{-1cm}\begin{CD}
(x^{(0)},x^{(1)},...,x^{(n-1)})&\longmapsto&(x^{(0)},x^{(1)},...,x^{(k)},0,...,0)&\longmapsto&(0,x^{(1)},x^{(2)},...,x^{(k)},0,...,0)\\
\mathbb{A}^n@>>>\mathbb{A}^n@>>>\mathbb{A}^n\\
@AAA@AAA@AAA\\
Y@>p_k>>\mathbb{A}^{k+1}@>q_k>>\mathbb{A}^k\\
&&(x^{(0)},x^{(1)},...,x^{(k)})&\longmapsto&(x^{(1)},x^{(2)},...,x^{(k)})
\end{CD}$$
Notons $\pi_k=q_k\circ p_k$. Pour établir ces résultats, Kronecker donne une méthode effective de décomposition d'une variété en composantes irréductibles, qui consiste à éliminer les inconnues du système d'équations successivement l'une après l'autre. Chaque composante peut enfin être représentée comme lieu des zéros d'un système d'au plus $(n+1)$ équations.
\\\par
Kronecker commence par décrire sa méthode d'élimination. Supposons $x^{(n-1)}$, $x^{(n-2)}$,...,$x^{(k+1)}$ déjà éliminées, et l'image par $p_k$ des composantes de $X$ de dimension $\leq k$ définie par un idéal $(H_1,H_2,...,H_m)$ qui, outre l'image de ces composantes, peut encore s'annuler sur certaines variétés <<~immergées~>>%
\footnote{\label{VarietesImmergees}Macaulay les appellera ainsi. Kronecker semble ne pas être conscient de l'existence de ces variétés immergées, dont Macaulay donnera un exemple (\cf \cite{macaulay1916}, exemple (\textit{ii}) p.~22).}, %
images de variétés incluses dans des composantes de $X$ de dimension $>k$. 
\par
Le choix d'un système de coordonnées suffisamment général permet d'assurer que l'image de toute composante irréductible de dimension $\leq k$ de $X$ par $p_k$ est un fermé, que $\pi_k$ fait de toute composante irréductible de $X$ de dimension $k$ un revêtement ramifié de $\mathbb{A}^k$, et que les images par $p_k$ de deux composantes de dimension $\leq k$ distinctes sont distinctes. Démontrons la première de ces trois assertions. La projection $p_k$ vérifie le diagramme commutatif suivant~:
$$\hspace{-2.8cm}
\begin{CD}
(Z:X^{(0)}:X^{(1)}:...:X^{(n-1)})&\longmapsto&\lbrace X^{(k+1)}=X^{(k+2)}=...=X^{(n-1)}=0\rbrace\cap L_{(Z:X^{(0)}:X^{(1)}:...:X^{(n-1)})}\\
\mathbb{P}^n-\lbrace Z=X^{(0)}=X^{(1)}=...=X^{(k)}=0\rbrace@>>>\mathbb{P}^n\\
@VVV@VVV\\
\mathbb{A}^n@>p_k>>\mathbb{A}^n\\
(\frac{X^{(0)}}{Z},\frac{X^{(1)}}{Z},...,\frac{X^{(n-1)}}{Z})&\longmapsto&(\frac{X^{(0)}}{Z},\frac{X^{(1)}}{Z},...,\frac{X^{(k)}}{Z},0,...,0)
\end{CD}$$
où $L_{(Z:X^{(0)}:X^{(1)}:...:X^{(n-1)})}$ est la sous-variété linéaire projective de codimension $(k+1)$ contenant $(Z:X^{(0)}:X^{(1)}:...:X^{(n-1)})$ et $\lbrace Z=X^{(0)}=X^{(1)}=...=X^{(k)}=0\rbrace$. Notons $\overline{Y}$ la clôture projective d'une composante $Y$ de $X$ de dimension $\leq k$. Alors 
$$\dim\overline{Y}\cap\lbrace Z=0\rbrace\leq k-1.$$ 
Or 
$$\dim\lbrace Z=X^{(0)}=X^{(1)}=...=X^{(k)}=0\rbrace=n-k-2$$
Donc, d'après le théorème de Bertini, pour un choix de coordonnées $X^{(0)}$,...,$X^{(k)}$ suffisamment générales, on aura~:
$$\dim\overline{Y}\cap\lbrace Z=X^{(0)}=X^{(1)}=...=X^{(k)}=0\rbrace<0$$
et $\overline{Y}$ évite $\lbrace Z=X^{(0)}=X^{(1)}=...=X^{(k)}=0\rbrace$. L'image de $Y$ par $p_k$ est donc fermée.
\par
On sépare de $(H_1,H_2,...,H_m)$ ses composantes de dimension $k$ en cherchant le plus grand commun diviseur des polynômes $H_i$, soit $F_{n-k}$ (ce qui suppose un algorithme de factorisation dans les anneaux de polynômes, dont on a déjà parlé). Pour tout $i$, notons aussi $K_i=\frac{H_i}{F_{n-k}}$. On obtient ainsi un idéal $(K_1,K_2,...,K_m)$ contenant l'image par $p_k$ des composantes de dimension $\leq (k-1)$. On élimine ensuite $x^{(k)}$. 
Pour ce faire, Kronecker calcule simplement le résultant des deux polynômes en $x^{(k)}$
$$\left\lbrace\begin{array}{l}
U_1K_1+U_2K_2+...+U_mK_m\\
V_1K_1+V_2K_2+...+V_mK_m
\end{array}\right.$$
où $U_1,...,U_m$ et $V_1,...,V_m$ sont des indéterminées. Ce résultant est une forme en ces indéterminées, dont les coefficients engendrent à leur tour un idéal dont le lieu des zéros contient l'image par $p_{k-1}$ des composantes de $X$ de dimension $\leq k-1$ (ainsi que, peut-être, certaines variétés <<~immergées~>>). Conformément aux idées de Kronecker,
pour étudier les composantes de dimension $k$, il faut <<~maintenir la variabilité~>> des grandeurs $x^{(1)}$, $x^{(2)}$,...,$x^{(k)}$. Le polynôme $F_{n-k}\in k(x^{(1)},...,x^{(k)})[x^{(0)}]$ fournit alors une équation
$$F_{n-k}=0$$
décrivant la variabilité de la grandeur $x^{(0)}$ dans les composantes de dimension $k$, et cette équation a un nombre fini de racines dans $\overline{k(x^{(1)},...,x^{(k)})}$. Autrement dit, l'hypersurface d'idéal $(F_{n-k})$ est conçue comme un revêtement ramifié de $\lbrace x^{(0)}=0\rbrace\subset\mathbb{A}^{k+1}$ (le choix de coordonnées suffisamment générales permet d'assurer que $x^{(0)}$ ne prend pas de valeur infinie).
\\\par
Pour étudier de manière effective la relation entre chaque composante irréductible $Y$ de $X$ de dimension $k$ et son image $p_k(Y)$, et en particulier la variabilité des autres variables $x^{(k+1)}$,...,$x^{(n-1)}$ en fonction de $x^{(1)}$,...,$x^{(k)}$, Kronecker utilise le même changement de variables que Poisson en 1802~; il choisit comme nouveau système de coordonnées
$$x=u_0x^{(0)}+u_1x^{(1)}+...+u_{n-1}x^{(n-1)},\ x^{(1)},\ x^{(2)},\ ...,\ x^{(n-1)}$$
où $u_0,u_1,...,u_{n-1}$ sont des indéterminées%
\footnote{Ici, nous changeons légèrement les notations de Kronecker. De plus, dans les \textit{Grundzüge}, Kronecker parle en fait de \textit{deux} changements de variables en même temps (nous faisions allusion au premier en disant qu'il fallait choisir un système de coordonnées <<~suffisamment général~>>). On pourrait croire qu'un deuxième changement de variables est superflu, puisque on a déjà supposé le premier suffisamment général pour éviter les cas d'exception. Mais le rôle du second changement de variables, comme l'a bien vu Macaulay, n'est pas d'éviter des cas d'exception. Il s'agit plutôt d'utiliser l'auxiliaire des coefficients indéterminés pour donner une expression algébrique de la relation entre $Y$ et $p_k(Y)$.
\Cf \cite{macaulay1916} p.~18, et aussi \cite{molk1992}, II.2, p.~154 et 158.}. %
On applique alors seulement la méthode d'élimination que nous venons de décrire, et on obtient
$$F_{n-k}\in k(u_0,u_1,...,u_{n-1},x^{(1)},x^{(2)},...,x^{(k)})[x]$$
$F_{n-k}$ est appelée <<~résolvente partielle~>>%
\footnote{Nous avons choisi, comme certains auteurs, l'orthographe française <<~résolvente~>> plutôt que <<~résolvante~>>, car le mot <<~résolvante~>> est utilisé dans un autre contexte pour désigner les <<~résolvantes de Lagrange~>> et les <<~résolvantes de Galois~>>.}. %
Kronecker appelle <<~résolvente totale~>> (\textit{Gesammtresolvente}) le produit $\prod F_i=0$ des résolventes partielles.
A chaque racine $\xi$ de $F_{n_k}$ correspond un feuillet de l'une des composantes de dimension $k$, paramétré par les variables $x^{(1)}$,...,$x^{(k)}$, et dont les points ont pour coordonnées $(\xi,x^{(1)},...,x^{(k)},\xi^{(k+1)},\xi^{(k+2)},...,\xi^{(n-1)})$. Ces points ont pour coordonnées, dans le système de coordonnées initial~:
$$(\xi^{(0)},x^{(1)},...,x^{(k)},\xi^{(k+1)},...,\xi^{(n-1)})$$
où l'on pose $\xi^{(0)}=\frac{1}{u_0}(\xi-u_1x^{(1)}-...-u_kx^{(k)}-u_{k+1}\xi^{(k+1)}-...-u_{n-1}\xi^{(n-1)})$. On peut donc écrire~:
$$\hspace{-1cm}F_{n-k}=\prod_i(x-\xi_i)=\prod_i\left(x-(u_0\xi_i^{(0)}+u_1x^{(1)}+...+u_kx^{(k)}+u_{k+1}\xi_i^{(k+1)}+...+u_{n-1}\xi_i^{(n-1)}\right)$$
où, pour tout $i$,
$$(\xi_i^{(0)},\xi_i^{(k+1)},...,\xi_i^{(n-1)})\in\left(\overline{k(u_0,...,u_{n-1},x^{(1)},...,x^{(k)})}\right)^{n-k}$$
On a ainsi justifié l'affirmation de Kronecker, que les $F_{n-k}$ sont <<~décomposables en facteurs linéaires~>>%
\footnote{\cf \cite{kronecker1882} p.~30.} %
en $u_0$,...,$u_{n-1}$. Le choix initial de coordonnées suffisamment générales permet encore d'assurer que les $\xi_i^{(0)}$,$\xi_i^{(k+1)}$,...,$\xi_i^{(n-1)}$ sont toutes finies, et l'on a ainsi décrit la réunion des composantes de dimension $k$ de $X$ comme revêtement ramifié de $\mathbb{A}^k$ par $\pi_k$. Mais il est possible que certains des feuillets
$$(\xi^{(0)},x^{(1)},...,x^{(k)},\xi^{(k+1)},...,\xi^{(n-1)})$$
ainsi décrits dépendent de $u_0,...,u_{n-1}$. \label{MacaulayLacuneKro}Macaulay a montré que si un tel feuillet existe, il appartient à une variété <<~immergée~>>. Donc toutes les composantes $Y$ de dimension $k$ de $X$ sont décrites par des paramétrages indépendants de $u_0,...,u_{n-1}$~:
$$(\xi_i^{(0)},\xi_i^{(k+1)},...,\xi_i^{(n-1)})\in\left(\overline{k(x^{(1)},...,x^{(k)})}\right)^{n-k}$$
\\\par
Un produit de $d$ facteurs $(x-\xi)$ décrivant tous les feuillets d'une même composante irréductible $Y$ est un facteur irréductible de $F_{n-k}\in k(u_0,...,u_{n-1})[x,x^{(1)},...,x^{(k)}]$. Substituons à nouveau à $x$, dans ce facteur, sa valeur
$$x=u_0x^{(0)}+u_1x^{(1)}+...+u_{n-1}x^{(n-1)},$$
on obtient alors un polynôme $\Phi\in\left(k[u_0,...,u_{n-1}]\right)[x^{(0)},x^{(1)},...,x^{(n-1)}]$. Kronecker écrit~:%
\footnote{\cf \cite{kronecker1882} p.~30. Kronecker ne démontre pas cette assertion. Selon Macaulay, König l'a démontrée dans \cite{konig1903} p.~234. Selon Netto et Le Vavasseur, \cite{molk1992} I~9~\S~69, K. Th. Vahlen aurait montré que $(n+1)$ équations sont parfois \textit{nécessaires}. Mais l'exemple qu'il donne est en fait erroné, \cf les cours de S. S. Abhyankar à Montréal (1970) sur le <<~nombre minimal d'équations définissant une courbe algébrique gauche~>>, \cite{abhyankar1971} p.~12.}
\begin{quote}
 En donnant aux grandeurs $u$ un certain nombre de système de valeurs, on voit facilement que $(n+1)$ tels systèmes de valeurs suffisent toujours à faire en sorte que les équations ainsi obtenues $\Phi_1=0$, $\Phi_2=0$,...,$\Phi_{n+1}=0$ aient pour résolvente totale $\Phi=0$.
\end{quote}
En effet%
\footnote{La méthode démonstrative exposée ici nous a été enseignée par Marc Giusti. Le choix des valeurs en lesquelles nous spécialisons les indéterminées peut être rendu effectif (grâce aux techniques exposées dans \cite{giusti1991a}).}, %
pour $1\leq i\leq n+1$, remplaçons les indéterminées $u_0,...,u_{n-1}$ respectivement par d'autres indéterminées $u_0^{(i)},...,u^{(i)}_{n-1}$, on obtient un polynôme $\Phi_i$ qui s'annule identiquement sur $Y$. On peut voir la variété d'idéal $(\Phi_1)$ comme un fibré sur l'espace affine de coordonnées $u^{(1)}_0,...,u^{(1)}_{n-1}$. Ecrivons $V(\Phi_1)$ comme réunion de ses composantes irréductibles\label{TrucGiusti}
$$V(\Phi_1)=\bigcup_jY^{(1)}_j\cup\bigcup_kZ^{(1)}_k$$
de sorte que les $Y^{(1)}_j$ soient les composantes irréductibles dont la fibre en $u^{(1)}_0,...,u^{(1)}_{n-1}$ est constante, et que les $Z^{(1)}_k$ soient les composantes irréductibles dont la fibre varie avec $u^{(1)}_0,...,u^{(1)}_{n-1}$. Ecrivons aussi $V(\Phi_1,\Phi_2)$ comme réunion de ses composantes irréductibles, en adoptant la même convention~:
$$V(\Phi_1,\Phi_2)=\bigcup_jY^{(2)}_j\cup\bigcup_kZ^{(2)}_k$$
Si l'on spécialise à présent $u^{(1)}_0,...,u^{(1)}_{n-1}$ en des valeurs quelconques dans $k$, il est possible de spécialiser $u^{(2)}_0,...,u^{(2)}_{n-1}$ dans $k$ de sorte que~:
$$\dim\bigcup_kZ^{(2)}_k=\dim\bigcup_kZ^{(1)}_k-1$$
Ecrivons de même, pour tout $1\leq i\leq n+1$,
$$V(\Phi_1,...,\Phi_i)=\bigcup_jY^{(i)}_j\cup\bigcup_kZ^{(i)}_k$$
On peut alors spécialiser les autres indéterminées $u^{(3)}_0,...,u^{(3)}_{n-1}$,...,$u^{(n+1)}_0,...,u^{(n+1)}_{n-1}$ de proche en proche, de sorte que~:
$$n>\dim\bigcup_kZ^{(1)}_k>\dim\bigcup_kZ^{(2)}_k>...>\dim\bigcup_kZ^{(n+1)}_k$$
Alors nécessairement $\bigcup_kZ^{(n+1)}_k=\varnothing$, et $V(\Phi_1,...,\Phi_{n+1})=Y$.
\\\par
Avec Macaulay, on peut aussi utiliser $\Phi$ pour décrire l'équivalence birationnelle entre $Y$ et $p_k(Y)$. Macaulay remarque%
\footnote{\cf \cite{macaulay1916} p.~27--28.} %
que le coefficient de $u_0^d$ dans $\Phi$ donne une équation de $p_k(Y)$~:
$$\phi(x^{(0)},x^{(1)},...,x^{(k)})=0$$
Pour tout $k+1\leq i\leq n-1$, le coefficient de $u_iu_0^{d-1}$ donne une équation de la forme~:
$$x^{(i)}\phi'-\phi_i=0,$$
où $\phi',\phi_i\in k[x^{(0)},x^{(1)},...,x^{(k)}]$. Cela donne un sens précis à l'assertion faite par Kronecker dans un langage encore vague à la fin de cette section~:%
\footnote{\textit{es lässt sich daher jedes $m$-fach ausgedehnte Gebilde einer beliebig grossen Mannigfaltigkeit auf ein solches eindeutig beziehen, das aus einer nur $(m+1)$-fachen Mannigfaltigkeit entnommen ist}, \cite{kronecker1882} p.~31.}
\begin{quote}
Pour toute composante de dimension $m$ d'une variété de dimension quelconque, il existe donc une correspondance biunivoque entre ses points et ceux d'une sous-variété d'une variété de dimension $(m+1)$.
\end{quote}

\section{Application~: la théorie de Galois}
Aux \S~11 et 12, Kronecker nous offre un premier exemple d'application de sa théorie. <<~La théorie de l'élimination, nous dit-il, montre la source propre de la nouvelle lumière qui fut apportée par Galois dans la théorie des équations algébriques, il y a un demi-siècle~>>. Le corps de décomposition d'une équation irréductible
$$ x^n-c_1x^{n-1}+...\pm c_n=0$$
à coefficients dans le domaine de base $\mathbb{Q}(\rat',\rat'',\rat''',...)$ est une extension finie, engendrée par ses racines $\xi_1, \xi_2,..., \xi_n$. L'on peut aussi interpréter ces racines comme la variété solution d'un système de $(n+i)$ équations, que Kronecker définit explicitement~:
$$(\bar{C}) \quad \fon_k(\xi_1,\xi_2,...,\xi_n)=c_k, \quad \mathfrak{g}_i(\xi_1,\xi_2,...,\xi_n)=c_i^\prime$$
où, pour $k=1,...,n$, les $\fon_k$ sont les fonctions symétriques élémentaires des racines, et $\mathfrak{g}_i$ sont des fonctions rationnelles bien choisies%
\footnote{\Cf \S~12. On remarquera que la variété des solutions est un ensemble fini de points, obtenus en faisant agir le groupe de l'équation sur les coordonnées de l'un quelconque d'entre eux. Mais vue comme variété sur $\mathbb{Q}(\rat',\rat'',\rat''',...)$, elle est irréductible. Par contre, la variété solution des équations $\fon_k(\xi_1,\xi_2,...,\xi_n)=c_k$, qui consiste en l'ensemble des points dont les coordonnées sont une permutation quelconque des racines, n'est en général pas irréductible sur $\mathbb{Q}(\rat',\rat'',\rat''',...)$. Les équations supplémentaires $\mathfrak{g}_i(\xi_1,\xi_2,...,\xi_n)=c_i^\prime$ sont donc nécessaires pour séparer les composantes irréductibles. Certainement Kronecker était-il conscient, en quelque manière, d'une telle interprétation géométrique (voir nos remarques plus loin sur ce sujet). Si l'on conçoit les coefficients de l'équation de départ comme étant des variables, on obtient un <<~revêtement galoisien~>>.}
. %
Conformément aux résultats du paragraphe précédent, Kronecker montrera plus loin (à la fin du \S~12) que $(n+1)$ équations suffisent en général. La résolvente du système d'équations est de la forme~:
$$g(x,u_1,u_2,...,u_n)=\prod (x-u_1\xi_{r_1}-u_2\xi_{r_2}-...-u_n\xi_{r_n})$$
où le produit est étendu à toutes les permutations $(r_1,...,r_n)$ appartenant au groupe de l'équation. En changeant les indices, on trouve~:
$$g(x,u_1,u_2,...,u_n)=\prod (x-u_{r_1}\xi_1-u_{r_2}\xi_2-...-u_{r_n}\xi_n)$$
et l'on voit ainsi que $g$ est une fonction des indéterminées $u_1, u_2,...,u_n$ \textit{invariante par les permutations du groupe}. Kronecker se flatte d'avoir ainsi donné à la théorie de Galois un <<~perfectionnement formel~>>, en ayant remplacé la considération du groupe par celle d'une <<~fonction concrète, invariable par le groupe~>>. Le lecteur moderne s'étonnera de cette préférence pour une <<~fonction concrète~>>, justement de la part de Kronecker, qui fit {\oe}uvre de pionnier dans l'axiomatisation de la notion de groupe abstrait%
\footnote{\Cf \cite{wussing1969}, pp.~44--48.}! %
Mais le calcul de la résolvente donne ici encore une méthode \textit{effective} pour déterminer le groupe de l'équation, étant donné la représentation de son corps de décomposition comme quotient d'un anneau de polynômes~:
$$\hspace{-1cm}\mathbb{Q}(\fon_1,...\fon_n)[x_1,...,x_n]/(\fon_1(x_1,...,x_n)-c_1,...,\fon_n(x_1,...,x_n)-c_n,\mathfrak{g}_1(x_1,...,x_n)-c_0^\prime,...,\mathfrak{g}_i(x_1,...,x_n)-c_i^\prime)$$
Bien sûr, tout le problème est alors de déterminer les $\mathfrak{g}_i$, mais Kronecker contourne d'abord ce problème en remarquant que la résolvente $g$ est un facteur irréductible de la résolvente du système d'équations
$$(C) \quad \fon_k(x_1,...,x_n)=c_k \quad (k=1,...,n)$$
qui s'écrit~:
$$\prod (x-u_1\xi_{r_1}-u_2\xi_{r_2}-...-u_n\xi_{r_n})$$
où le produit est cette fois étendu à toutes les permutations $(r_1,...,r_n)$ du groupe symétrique. La théorie de Galois est donc ramenée à un problème de factorisation. Enfin les coefficients de $g$ sont des fonctions des racines invariantes par les permutations du groupe, et ils fournissent d'ailleurs le système d'équations $\mathfrak{g}_i(\xi_1,...,\xi_n)=c_i^\prime$ cherché. Mais comme dans la fonction $g$, les rôles des indéterminées $u$ et des racines $\xi$ sont interchangeables, Kronecker préfère se ramener à l'étude des fonctions rationnelles des $n$ indéterminées $u$ dans le \S~12~: après adjonction des indéterminées $u_1,...,u_n$, le corps de décomposition de l'équation de départ est devenu corps de rupture de l'équation $g(x)=0$ sur $\mathbb{Q}(\rat',\rat'',\rat''',...,u_1,u_2,...,u_n)$. Ce dernier corps $\mathbb{Q}(\rat',\rat'',\rat''',...,u_1,u_2,...,u_n)$ est lui-même corps de décomposition d'une équation en $u$, sur un sous-corps contenant les fonctions symétriques élémentaires de $u_1,...,u_n$ et les coefficients de $g$. Le groupe de Galois de ce nouveau corps de décomposition est le même que le groupe de Galois de l'équation de départ. On a donc ramené la question à l'étude du corps des fonctions rationnelles de $n$ indéterminées $u_1,...,u_n$.  
Le paragraphe \S~11 se présente donc comme une incursion rapide, et difficile à lire, en théorie de Galois~: Kronecker n'écrit jamais le système d'équations $(\bar{C})$ (bien qu'il utilise cette notation pour le désigner), et une erreur typographique change le $(\bar{C})$ en $(C)$ (fin de la page 33).
\par
Après avoir ainsi ramené la question à l'étude du corps des fonctions rationnelles de $n$ indéterminées, Kronocker propose (au \S~12) de prendre à présent <<~comme point de départ~>>, en théorie de Galois, l'étude de ce corps. Cela revient justement à faire la théorie de Galois d'une équation \textit{générique} de degré $n$. Le corps des fonctions rationnelles à $n$ indéterminées sur $\mathbb{Q}$ peut en effet se représenter comme extension finie du corps des fonctions symétriques, comme le montre le diagramme suivant, où l'on a aussi indiqué les anneaux d'entiers étudiés par Kronecker~:
\begin{equation*}
\begin{CD}
\begin{array}{l@{}l}
K&=\mathbb{Q}(\fon_1,...,\fon_n)[x_1,...,x_n]/(\fon_1-\sum x_i,...,\fon_n-\prod x_i)\\
&=\mathbb{Q}(x_1,...,x_n)
\end{array} @<<< O_K=\mathbb{Q}[x_1,...,x_n]\\
@AAA   @AAA\\
k=\mathbb{Q}(\fon_1,...,\fon_n) @<<< O_k=\mathbb{Q}[\fon_1,...,\fon_n]
\end{CD}
\end{equation*}
$K$ est de degré $n!$ sur $k$. Kronecker montre que le $O_k$-module $O_K$ est libre de rang $n!$, en exhibant une base%
\footnote{On a donc un nouveau cas particulier dans lequel le nombre de générateurs d'une base de $O_K$ est égal au degré de $K$ sur $k$ (Kronecker l'avait déjà annoncé au \S~7). \note{Mais il ne démontre pas que l'anneau $O_K$ est la clôture intégrale de $O_k$ dans $K$.}}
$$\left\lbrace x_1^{h_1}...x_{n-1}^{h_{n-1}}\right\rbrace_{\begin{array}{l} h_k=0,...,n-k\\k=1,...,n-1 \end{array}}$$
Le \textit{discriminant du genre} est alors le carré du déterminant
$$\left\vert\left(x_{\sigma(1)}^{h_1}...x_{\sigma(n-1)}^{h_{n-1}}\right)_{h_1,...,h_{n-1},\sigma}\right\vert$$
où $h_k=0,...,n-k$ et $\sigma$ parcourt le groupe symétrique $\Sigma_n$. En ordonnant convenablement les éléments de la base, on peut réécrire cette matrice comme produit $AB$ des deux matrices suivantes. La matrice $A$ est une matrice diagonale par blocs, le $i$-ème bloc sur la diagonale étant de la forme
$$\left(x_{\sigma\tau^i(2)}^{h_2}...x_{\sigma\tau^i(n-1)}^{h_{n-1}}\right)_{h_2,...,h_{n-1},\sigma}$$
où $h_k=0,...,n-k$, où $\tau$ est la permutation cyclique $(12...n)$, et où $\sigma$ parcourt le sous-groupe de $\Sigma_n$ laissant $i$ invariant. La matrice $B$ est constituée elle aussi de $n\times n$ blocs de taille $(n-1)!\times (n-1)!$, le bloc $(i,j)$ étant simplement $x_j^{i-1}\text{Id}_{(n-1)!}$. Notons $\mathcal{D}$ le carré du déterminant de Vandermonde (c'est-à-dire le discriminant de l'équation). On vérifie alors aisément que $\det B=\pm\mathcal{D}^{\frac{(n-1)!}{2}}$. Par récurrence sur $n$, on montre alors que
$$\left\vert\left(x_{\sigma(1)}^{h_1}...x_{\sigma(n-1)}^{h_{n-1}}\right)_{h_1,...,h_{n-1},\sigma}\right\vert^2=\mathcal{D}^{\frac{n!}{2}}$$
comme l'annonce Kronecker sans le démontrer à la fin du $\S~12$.
\par
De même, si $L=\mathbb{Q}(\fon_1,...,\fon_n,\mathfrak{g})$ est une sous-extension de degré $\rho$ sur $k$, Kronecker démontre que $\rho$ générateurs suffisent à engendrer le $O_k$-module $O_L=L\cap O_K$. Si l'équation
$$x^n-\fon_1x^{n-1}+...\pm\fon_n=0$$
reste irréductible sur $L$, Kronecker dit que l'extension $L$ est <<~propre~>> (\textit{eigentlich}). Dans ce cas, le groupe de l'équation sur $L$ est d'ordre $r=\frac{n!}{\rho}$, et Kronecker propose, plutôt que de parler du <<~groupe de l'équation~>> (notion qui renvoie à la nature de l'indétermination des racines), d'introduire le terme <<~genre d'affection~>> (\textit{Affect-Gattung}), pour désigner le corps des fonctions rationnelles de $n$ indéterminées invariantes par les permutations du groupe, notion qui renvoie donc, en quelque sorte, à la nature arithmétique des racines et des fonctions rationnelles des racines. L'équation est alors dite avoir une certaine <<~affection~>> d'ordre $r$, et ses racines appartenir à une certaine <<~classe~>> d'ordre $r$ (\cf \cite{kronecker1882}, p.~37). Ce faisant, Kronecker se réclame plus d'Abel que de Galois, et mentionne à cet égard un problème de Galois inverse~:
\begin{quote}
  Il échappe à Galois un problème des plus intéressants, dans la théorie des équations algébriques, qu'Abel a trouvé et aussi traité. C'est le problème de déterminer toutes les équations d'une certaine classe pour un domaine de rationalité donné, et je vais ici l'exposer plus en détail, aussi parce qu'il montre clairement la nature arithmétique des questions algébriques.%
\footnote{\Cf \cite{kronecker1882}, p.~38.}
\end{quote}
Si $\mathfrak{g}$ est un élément primitif de $L$ sur $k$, déterminer les équations de genre d'affection $L$ (c'est-à-dire les équations dont le groupe de Galois est le groupe de $K$ sur $L$), à coefficients dans un corps donné, revient à trouver certaines des solutions dans ce corps d'une équation $\Phi(\mathfrak{g},\fon_1,...,\fon_n)=0$ (qui n'est autre que l'équation de $\mathfrak{g}$ sur $k$). Il s'agit donc, remarque Kronecker, d'un <<~problème diophantien~>>. Pour mieux comprendre cette attitude, il suffira de rappeler ses travaux sur le théorème dit <<~de Kronecker-Weber~>>. Ces travaux s'inscrivent dans la continuité des recherches menées par Abel dans son manuscrit \textit{Sur la résolution algébrique des équations}.
Abel y proposait un certain problème de Galois inverse~: <<~Trouver toutes les équations d'un degré déterminé quelconque qui soient résolubles algébriquement~>>%
\footnote{\cf \cite{abel1881} t.~2 p.~219. Voir aussi \cite{houzel2004} p.~62--65\note{Epreuve p.~33}.}. %
Le théorème de Kronecker-Weber résout le problème de Galois inverse suivant~: quelle est la forme générale d'une extension galoisienne de $\mathbb{Q}$ dont le groupe est abélien~? En 1853, Kronecker affirmait déjà que toute extension de $\mathbb{Q}$ dont le groupe est cyclique se plonge dans une extension engendrée par des racines de l'unité%
\footnote{Kronecker écrit dans \cite{kronecker1853} p.~10 que <<~la racine de toute équation abélienne à coefficients entiers peut être représentée comme fonction rationnelle de racines de l'unité~>>. Mais ce sont seulement les équations de groupe cyclique qu'il désigne alors par l'expression <<~équation abélienne~>>.}, %
et il en donnait l'esquisse d'une démonstration (Weber en donnerait plus tard une démonstration complète). En 1877, il généralise et énonce le théorème qui a gardé son nom, et souligne son intérêt pour la théorie des nombres~:
\begin{quote}
  Toutes les racines des équations abéliennes à coefficients entiers sont fonctions rationnelles de racines de l'unité, et toutes les fonctions rationnelles de racines de l'unité sont racines d'équations abéliennes à coefficients entiers.\\
Cette proposition donne, me semble-t-il, un précieux aperçu en théorie des nombres algébriques~; car elle constitue un premier progrès quant à leur classification naturelle, qui dépasse leur rassemblement en genres, jusque aujourd'hui le seul que l'on ait considéré.%
\footnote{\Cf \cite{kronecker1877}, p.~69. Dans cet article, Kronecker utilise l'expression \textit{Abelsche Gleichung} pour désigner une équation abélienne, et l'expression \textit{einfache Abelsche Gleichung} pour désigner une équation dont le groupe est cyclique.}
\end{quote}
Dans ce même article, il affirme qu'il avait été conduit, en 1857, dans ses recherches sur les fonctions elliptiques et la multiplication complexe, à formuler une conjecture analogue sur les extensions abéliennes d'un corps quadratique imaginaire. Nous reviendrons sur ce sujet, étroitement lié à la théorie du corps de classes, quand nous commenterons la section \S~19 des \textit{Grundzüge}.

\section{Théorie de la divisibilité dans les anneaux d'entiers algébriques}
On a vu ci-dessus, dans la construction de la \textit{résolvente} d'un système d'équations par Kronecker, une première utilisation de l'<<~auxiliaire des coefficients indéterminés~>> (\textit{methodische Hülfsmittel der unbestimmten Coefficienten}). Il utilise à nouveau cette méthode de l'adjonction d'indéterminées au domaine de base dans les sections \S~14 à 18, pour construire une théorie de la divisibilité dans les anneaux d'entiers. Pour qualifier cette méthode, il parlera  plus loin d'un accroissement de la <<~dimension~>>, en utilisant un langage géométrique~:
\begin{quote}
De même que la ligne des nombres réels s'étend par une <<~unité latérale~>> en le plan des nombres complexes, un domaine de grandeurs $[\rat',\rat'',\rat''',...]$ est agrandi par les indéterminées [...], d'une certaine manière, en rapport à sa <<~dimension~>>.%
\footnote{\Cf \S~22. \note{Voir ici les autres passages des \textit{Grundzüge} qui mentionne une interprétation géométrique.}}
\end{quote}
Il ne s'agit donc pas seulement d'un auxiliaire algébrique formel. Cette méthode est à la fois censée être la plus simple, mais aussi produire la description la plus simple de la réalité mathématique. Kronecker la compare à l'un des grands progrès des mathématiques du XIX\up{ème} siècle~: bien des phénomènes de l'analyse gagnent en <<~simplicité~>> quand ils sont étudiés du point de vue de l'analyse complexe%
\footnote{\Cf \cite{kronecker1882}, p.~94.}. %
Quelle que soit l'interprétation suggérée ici par Kronecker, ses successeurs ont d'ailleurs naturellement adopté une interprétation géométrique des \textit{Grundzüge}, en particulier de sa méthode d'élimination, comme on l'a déjà remarqué. Donc, bien que le langage de Kronecker soit rarement géométrique, le plus souvent arithmético-algébrique, il semble légitime d'en pousser l'interprétation au delà du formel, bien qu'il faille avoir recours pour cela au substrat intuitif auquel nous sommes aujourd'hui habitués, celui de la théorie des ensembles et des théories de structures (telle la notion d'idéal). C'est ce que l'on a fait ci-dessus en traduisant son langage dans les notations modernes de l'algèbre commutative. Pour faire une théorie de la divisibilité dans les anneaux d'entiers, qui ne sont pas toujours factoriels, l'on a recours aujourd'hui à la théorie des idéaux de Dedekind. Dans les \S~14 à 18, Kronecker propose une autre solution, commentée de manière détaillée par Harold M. Edwards dans \cite{edwards1980}. Nous allons seulement résumer son analyse, avant de proposer notre propre interprétation%
\footnote{H. Weyl a proposé encore une autre lecture, axiomatique, de la théorie des diviseurs de Kronecker, dans \cite{weyl1940}, chapitre 2. Il dresse une liste d'axiomes que doit vérifier le concept de <<~diviseur~>>. Puis il définit le diviseur $\mathfrak{A}$ comme une suite finie d'entiers algébriques $(\alpha_1,...,\alpha_r)$, en donnant un critère de divisibilité d'un entier algébrique quelconque $\alpha$ par $\mathfrak{A}$, et il montre que cette définition vérifie les axiomes (p.~49).}. %
Kummer, dès 1844, avait construit une théorie de la divisibilité pour les entiers cyclotomiques (éléments de $\mathbb{Z}[\zeta_p]$, où $\zeta_p$ est une racine $p$-ième de l'unité, $p$ premier), en introduisant des <<~nombres idéaux~>>, afin de rétablir l'existence et l'unicité de la décomposition d'un entier en facteurs premiers. Un tel <<~nombre idéal~>> n'était pas défini autrement par Kummer qu'en énonçant la relation d'équivalence entre deux entiers cyclotomiques modulo ce nombre. Dedekind avait au contraire choisi comme fondement, à partir de 1871, l'<<~idéal~>>, ensemble des entiers divisibles par un <<~nombre idéal~>> donné. Il introduisit plus tard (années 1890) des idéaux fractionnaires, et donna ainsi à la théorie la forme qu'elle a gardée jusqu'à aujourd'hui. Kronecker choisit de définir, plutôt que des <<~nombres idéaux~>> ou des <<~idéaux~>>, des <<~diviseurs~>>, définis (sur le modèle de Kummer) par une relation d'équivalence. Seulement (au contraire de Kummer) il exprime cette relation d'équivalence comme une véritable relation de divisibilité, à condition d'adjoindre à l'anneau d'entiers des indéterminées. Un diviseur est noté $$\mod[\phi x+\phi' x'+\phi'' x''+...]$$ où $\phi x+\phi' x'+...$ est une forme à coefficients $x,x',...$ dans l'anneau d'entiers. Kronecker définit la norme $\Nm(\phi x+\phi' x'+...)$ d'une forme (produit des conjugués, la conjugaison agissant sur les coefficients de la forme). La forme est dite <<~primitive~>> si le plus grand commun diviseur des coefficients de sa norme est 1 (il faut que l'anneau d'entiers du corps de base soit factoriel pour que cela ait un sens). Kronecker définit aussi la forme primitive $\Fm(\phi x+\phi' x'+...)$, telle que~:
$$P.\Fm(\phi x+\phi' x'+...)=\Nm(\phi x+\phi' x'+...)$$
où $P$ est le plus grand commun diviseur des coefficients de la norme (appelé aujourd'hui son <<~contenu~>>). Un entier quelconque $z$ est alors dit <<~divisible~>> (\textit{theilbar}) par ce diviseur si et seulement si le quotient
$$\frac{z.\Fm(\phi x+\phi' x+...)}{\phi x+\phi' x+...}$$
est une forme $Q$ à coefficients entiers algébriques. H. M. Edwards a remarqué que les lecteurs de Kronecker (Dedekind, Hurwitz) ont buté sur ce point, Kronecker ne laissant pas clairement entendre si le quotient $Q$ doit être une forme à coefficients entiers algébriques, ou si $Q$ doit simplement vérifier une équation $Nm(X-Q)=0$ à coefficients entiers dans le corps de base (les deux étant équivalents, bien que Kronecker ne le démontre pas)%
\footnote{\Cf \cite{edwards1980}, p.~356 et 364--368. Ce problème amena Dedekind à concevoir son <<~théorème de Prague~>>, qui permet, comme le montre Edwards, de simplifier considérablement la théorie de Kronecker. \Cf \textit{infra} p.~\pageref{ChainonManquant}.}. %
Deux diviseurs sont dits <<~absolument équivalents~>> lorsqu'ils définissent la même relation de divisibilité. Kronecker construit ainsi toute une théorie de la divisibilité. A la section \S~14, il montre que le diviseur $\mod[x+u'x'+u''x''+...]$, défini par une forme \textit{linéaire}, est en un certain sens le plus grand commun diviseur des entiers $x,x',x'',...$. A la section \S~15, il définit le vocabulaire et les notations, puis il énonce et il démontre son <<~premier théorème fondamental~>>~:
\begin{quote}
Quand le produit de deux formes algébriques entières [=formes dont les coefficients sont des grandeurs algébriques entières], dont l'une est primitive, est congruent à zéro modulo un diviseur algébrique, alors l'autre forme doit elle-même être divisible par le diviseur.%
\footnote{\Cf \S~14, IX.}
\end{quote}
Les sections \S~16 et 17 sont consacrées à la démonstration d'un <<~second théorème fondamental~>>~:
\begin{quote}
Des diviseurs algébriques qui ont les mêmes éléments [=les mêmes coefficients] sont absolument équivalents.%
\footnote{\Cf \S~16, II.}
\end{quote}
En particulier tout diviseur est absolument équivalent à un diviseur défini par une forme \textit{linéaire}. Mais pour démontrer le <<~second théorème fondamental~>>, Kronecker commence justement par démontrer (\S~16) plusieurs propriétés des diviseurs définis par des formes linéaires, dont un analogue du lemme de Gauss~:
\begin{quote}
Si le produit de deux diviseurs algébriques est divisible par un troisième, et que le premier diviseur n'a aucun diviseur commun avec le troisième, alors le deuxième est divisible par le troisième.%
\footnote{\Cf \S~16, V.}
\end{quote}
H. M. Edwards remarque que Kronecker utilise, dans la démonstration du second théorème fondamental, l'existence, pour tout entier du corps de base, d'un diviseur qui le divise et soit défini par une forme linéaire, fait essentiel qu'il ne démontre proprement qu'à la section \S~18. Enfin, à la section \S~18, Kronecker décrit une méthode effective pour obtenir tous les diviseurs d'un entier du corps de base, et il en déduit l'existence et l'unicité de la décomposition en facteurs premiers d'un quelconque diviseur. 
\par
De ce point de vue, la théorie de Kronecker est donc une modification et une généralisation de celle de Kummer, en donnant le rôle central à la notion de plus grand commun diviseur (incarnée par les diviseurs définis par des formes linéaires, et qui apparaît dès le \S~14) plutôt qu'à celle de facteur premier (l'existence et l'unicité de la décomposition en facteurs premiers n'étant décrite qu'au \S~18, sans même qu'en soit proprement énoncé un théorème). Kronecker se flatte cependant d'avoir donné une existence <<~véritable~>> (\textit{wirklich}) aux diviseurs en les représentant par des formes, grâce à l'<<~auxiliaire des coefficients indéterminés~>>. L'on pourrait en donner l'interprétation moderne suivante.
\par
Soit une extension galoisienne $K$ sur $k$, avec des anneaux d'entiers $O_k$ et $O_K$~:
\begin{equation*}
\begin{CD}
K @<<< O_K\\
@AAA   @AAA\\
k @<<< O_k
\end{CD}
\end{equation*}
On suppose $O_k$ factoriel.
Soit $S\subset O_k[u_1,u_2,...]$ l'ensemble multiplicatif des <<~formes primitives~>> en les indéterminées $u_1,u_2,...$, à coefficients dans $O_k$ (les indéterminées forment un ensemble dénombrable $\lbrace u_1,u_2,... \rbrace$~; on emploiera dans tout ce qui suit le mot <<~forme~>> pour désigner un polynôme en $u_1,u_2,...$ \textit{non nécessairement homogène}, suivant ainsi l'habitude de Kronecker). Alors $S^{-1}O_k[u_1,u_2,...]$ est intégralement clos dans son corps de fractions $k(u_1,u_2,...)$, et $O_K\otimes S^{-1}O_k[u_1,u_2,...]$ est entier sur $S^{-1}O_k[u_1,u_2,...]$ \note{indiquer une référence pour ces deux résultats}. On a donc plongé $O_k$ et $O_K$ dans deux nouveaux anneaux d'entiers, et l'on a à présent la situation suivante~:
\begin{equation*}
\begin{CD}
K(u_1,u_2,...) @<<< O_K\otimes S^{-1}O_k[u_1,u_2,...]\\
@AAA   @AAA\\
k(u_1,u_2,...) @<<< S^{-1}O_k[u_1,u_2,...]
\end{CD}
\end{equation*}
\begin{theo}
\label{factoriel}
L'anneau $O_K\otimes S^{-1}O_k[u_1,u_2,...]$ est factoriel.
\end{theo}
\begin{corol}
Cet anneau est noethérien et principal%
\footnote{En fait, suivant Weyl (\cite{weyl1940}, p.~60), on peut même montrer que tout idéal est de la forme $(x_1u_1+x_2u_2)$ où $x_1, x_2 \in O_K$.}. %
\end{corol}
\textit{Démonstration}. On démontrera d'abord que tout idéal $(x_1,...,x_k)$ engendré par un nombre fini de générateurs est principal, engendré par $\text{pgcd}(x_1,...,x_k)=x_1u_1+...+x_ku_k$, où $u_1,...,u_k$ sont des indéterminées qui ne figurent pas dans $x_1,...,x_k$.
\par
On prolonge les morphismes de conjugaison de $K$ à $K(u_1,u_2,...)$, ce qui permet de définir une norme dans $K(u_1,u_2,...)$ (il n'y a même pas besoin que l'extension $K$ sur $k$ soit galoisienne pour cela). On définit de même $\Fm\alpha$, pour tout $\alpha\in O_K\otimes S^{-1}O_k[u_1,u_2,...]$, tel que~:
$$\Nm\alpha=P_\alpha.\Fm\alpha$$
où $P_\alpha$ est le plus grand commun diviseur des coefficients du numérateur de $\Nm\alpha$. Alors $\Fm\alpha$ est un quotient de formes primitives, donc une unité.
\par
Soit $\alpha=x_1u_1+...+x_ku_k$. Montrons que pour tout $i$, $\frac{x_i\Fm\alpha}{\alpha}$ est entier sur $S^{-1}O_k[u_1,u_2,...]$. Il suffit pour cela de voir que l'équation
$$\Nm\left(X-\frac{x_i\Fm\alpha}{\alpha}\right)=0$$
est une équation à coefficients entiers. Or,
$$\Nm\left(X-\frac{x_i\Fm\alpha}{\alpha}\right)=\frac{\Nm(X\alpha-x_i\Fm\alpha)}{P_\alpha.\Fm\alpha}$$
et l'on voit facilement que $\Nm(X\alpha-x_i\Fm\alpha)$ est bien divisible par $P_\alpha$ (\cf la section \S~14 de Kronecker). Donc pour tout $i$, $x_i$ est divisible par $\alpha$. Donc en termes d'idéaux, $(x_1,...,x_k)=(\alpha)$.
\par
Pour la suite de la démonstration, on aura besoin de deux lemmes. On dira qu'un diviseur de la forme $x_1u_1+...+x_ku_k$, avec $x_i\in O_K$ pour tout $i$, est un <<~diviseur linéaire~>>.
\begin{lemme}
\label{existdiv}
Soit $p$ premier dans $O_k$. Alors $p$ possède un diviseur linéaire premier qui n'est pas une unité.
\end{lemme}
\textit{Démonstration} (\cf section \S~18). Il existe au moins un diviseur linéaire de $p$ qui n'est pas une unité, à savoir $pu_0$. S'il en existe d'autres, on peut les écrire sous la forme $pu_0+x_1u_1+...+x_ku_k$ (à une unité près), car si $(x_1,...,x_k)$ est un diviseur de $p$ alors $(x_1,...,x_k)=(p,x_1,...,x_k)$. Soit $\Sigma$ l'ensemble de ces diviseurs. Tout diviseur linéaire de $p$ est alors, à une unité près, un élément de $\Sigma$. Si $O_k=\mathbb{Z}$, quitte à réduire modulo $p$ les entiers $x_1,..,x_k$, on voit que l'on peut d'ailleurs se limiter à un ensemble $\Sigma$ fini%
\footnote{Remarquons en effet que pour $O_k=\mathbb{Z}$, l'anneau des restes $O_K/(p)$ est de cardinal fini $p^{[K:k]}$.}. %
On ordonne $\Sigma$ par la relation de divisibilité. Montrons qu'il a un élément minimal. Soient $\alpha,\beta,\gamma\in\Sigma$, avec $\beta=\alpha\gamma$. Alors $P_\beta=P_\alpha P_\gamma$. Donc si $\alpha$ divise strictement $\beta$, $P_\gamma\neq 1$, et $P_\alpha$ divisera strictement $P_\beta$. Si $P_\alpha=1$, $\alpha$ est minimal, sinon  on recommence avec $\alpha$ ce que l'on vient de faire avec $\beta$. Comme $O_k$ est factoriel, on obtient ainsi un élément minimal $\alpha$ après un nombre fini d'étapes.%
\footnote{Remarquons toutefois que, dans le cas où $\Sigma$ est infini, cette méthode n'est pas effective car elle ne précise pas comment tester si $\alpha$ est minimal, ni comment trouver, sinon, un nouveau diviseur $\beta$. Kronecker annonce qu'il donnera une méthode meilleure, dans la section \S~25, pour la décomposition en facteurs premiers~; \note{comme on le verra} cette autre méthode marche dans le cas $O_k=\mathbb{Z}$, mais il semble douteux qu'elle soit applicable en général.}. %
Montrons à présent que $\alpha$ est premier. Soient $\beta,\gamma\in O_K\otimes S^{-1}O_k[u_1,u_2,...]$, tels que $\alpha\nmid\beta$, $\alpha\nmid\gamma$. Montrons que $\alpha\nmid\beta\gamma$. On fait une récurrence sur le nombre $m$ d'indéterminées présentes dans $\beta$ ou $\gamma$. On peut supposer qu'aucune de ces indéterminées n'apparaît dans $\alpha$, quitte à multiplier par une unité%
\footnote{Deux diviseurs linéaires ayant les mêmes coefficients sont égaux à une unité près. Il s'agit là d'un théorème énoncé par Kronecker au \S~16. En effet, pour tout $i$, $x_1v_1+...+x_kv_k$ divise $x_i$. Alors $x_1v_1+...+x_kv_k$ divise $x_1u_1+...+x_ku_k$.}. %
Si $m=0$, $\text{pgcd}(\beta,\alpha)=\beta u_1+\alpha u_2$ est un diviseur linéaire de $\alpha$, et, à cause du caractère minimal de $\alpha$, est donc trivial. Donc $\beta\gamma u_1+\alpha\gamma u_2=\epsilon\gamma$ où $\epsilon$ est une unité, et comme $\alpha$ ne divise pas $\gamma$, il ne divise pas non plus $\beta\gamma$. Si $m>0$, soit $u$ l'une des indéterminées présentes dans $\beta$ ou $\gamma$. On écrit~:
$$\beta=\sum\beta_iu^i, \quad \gamma=\sum\gamma_ju^j$$
Alors si $s$ et $t$ sont les plus grands indices tels que $\alpha\nmid\beta_s$ et $\alpha\nmid\gamma_t$, le coefficient de $u^{s+t}$ dans $\beta\gamma$ ne sera pas divisible par $\alpha$ (par récurrence sur $m$). Donc $\beta\gamma$ non plus.
\begin{lemme}[Lemme de Gauss]
\label{gauss}
Soient $p, a, b \in O_K\otimes S^{-1}O_k[u_1,u_2,...]$ avec $p$ premier, $p \nmid a$ et $p \mid ab$. Alors $p \mid b$.
\end{lemme}
\textit{Démonstration}. Par hypothèse, $\text{pgcd}(p,a)=pu_1+au_2=\epsilon$ est une unité. Donc $bpu_1+bau_2=b\epsilon$, et $p \mid ab \Longrightarrow p \mid b$.  
\begin{corol}
Aux unités près, il y a unicité de la décomposition en facteurs premiers dans $O_K\otimes S^{-1}O_k[u_1,u_2,...]$.
\end{corol}
\textit{Fin de la démonstration du théorème \ref{factoriel}} (existence de la décomposition en facteurs premiers). Soit $\alpha\in O_K\otimes S^{-1}O_k[u_1,u_2,...]$. Si $P_\alpha$ est une unité de $O_k$, $\alpha$ est aussi une unité. Sinon, soit $p\in O_k$ premier tel que $p \mid P_\alpha$. Il existe (lemme \ref{existdiv}) un diviseur linéaire premier $\pi$ de $p$, donc de $\Nm\alpha$. Alors (lemme \ref{gauss}) $\pi$ divise l'un des conjugués $\sigma(\alpha)$ de $\alpha$. Donc $\sigma^{-1}(\pi)$ est un diviseur linéaire premier de $\alpha$ (il faut ici utiliser l'hypothèse que $K$ est galoisien sur $k$~; \cf \textit{remarque} ci-dessous). Soit 
$$\beta=\frac{\alpha}{\sigma^{-1}(\pi)}$$
On a alors~:
$$P_\beta.P_{\sigma^{-1}(\pi)}=P_\alpha$$
Si $P_\beta=1$, c'est fini~; sinon l'on recommence avec $\beta$ ce que l'on a fait avec $\alpha$ et l'on conclut par récurrence sur le nombre de facteurs premiers de $P_\alpha$ (car $O_k$ est factoriel). \textit{q.~e.~d.}\\
\textit{Remarque~:} on peut sûrement se passer de l'hypothèse que l'extension $K$ est galoisienne dans la démonstration du théorème \ref{factoriel}, en plongeant $K$ dans une extension galoisienne.
\par
Dans la section \S~18, Kronecker annonce un théorème sur la ramification (il renvoie au \S~25 pour plus de détails). Un entier premier $p\in O_k$ est dit \textit{ramifié} si ses facteurs premiers dans \mbox{$O_K\otimes S^{-1}O_k[u_1,u_2,...]$} ne sont pas tous distincts. 
\begin{theo}
\label{ramif}
Un entier premier $p\in O_k$ est ramifié si et seulement s'il divise le discriminant du genre.
\end{theo}
A la suite de ce théorème, Kronecker énonce quelques assertions qui sont en contradiction avec sa définition de la norme, comme produit des conjugués \textit{distincts}. Désignons-donc par <<~norme absolue~>> le produit des $n$ conjugués, où $n=[K:k]$. Soit $p\in O_k$ premier. Sa norme absolue est $p^n$. Alors la norme absolue de chaque diviseur premier de $p$ divise $p^n$, et est donc elle-même de la forme $p^f$. L'exposant $f$ est appelé <<~ordre~>> du diviseur. Si $p$ n'est pas ramifié, la somme des ordres des diviseurs premiers de $p$ est donc égale à $n$. Lorsque $O_k=\mathbb{Z}$, le discriminant du genre est différent de 1. On a donc le théorème%
\footnote{Kronecker énonce ce théorème à la section 8, p.~21. C'est une conséquence immédiate de la méthode de factorisation qu'il expose dans la dernière section des \textit{Grundzüge}. \Cf \cite{edwards1990}, \S~2.8, qui attribue ce théorème à Dedekind. Abhyankar \cite{abhyankar1976} p.~444 l'attribue à Kronecker mais fait une remarque intéressante. Cassels, Fröhlich \cite{cassels1967} l'attribuent à Minkowski qui en a aussi donné une\note{autre ?} démonstration\note{géométrique ?}.}~:
\begin{theo}
Il n'existe pas d'extension finie non ramifiée de $\mathbb{Q}$.
\end{theo}
Sous les mêmes hypothèses, Kronecker montre que le nombres des restes dans $O_K$ modulo un diviseur est égal à la norme du diviseur. A la section 19, il continue son incursion en théorie des nombres et propose d'appliquer le <<~principe d'équivalence de Kummer~>>, c'est-à-dire de considérer, en termes modernes, le groupe des classes d'idéaux. Un diviseur, élément de $O_K\otimes S^{-1}O_k[u_1,u_2,...]$, est dit appartenir à la <<~classe principale~>> (\textit{Hauptclasse}) s'il appartient (à une unité près) à $O_K$. Le groupe des classes est alors défini comme quotient du monoïde multiplicatif de l'anneau $O_K\otimes S^{-1}O_k[u_1,u_2,...]$ par la classe principale. La relation d'équivalence sous-jacente entre diviseurs est appelée par Kronecker <<~équivalence relative~>>, pour la distinguer de la relation d'<<~équivalence absolue~>> entre deux diviseurs qui ne diffèrent qu'à une unité près. Kronecker donne deux critères équivalents d'équivalence relative. Notons $H$ la classe principale, et supposons l'extension galoisienne (hypothèse que Kronecker omet mais qui semble ici nécessaire pour garantir que le monoïde quotient soit un groupe~: l'inverse d'un diviseur, dans le groupe des classes, est en effet représenté par le produit de ses conjugués\note{Si l'on définit le groupe des classes au moyen des idéaux fractionnaires, y a-t-il au moins un idéal entier dans chaque classe?}). \\
\textit{1\up{er} critère d'équivalence relative}~: les diviseurs $\phi$ et $\psi$ sont équivalents si et seulement s'il existe un diviseur $\chi$ tel que $\phi\chi$ et $\psi\chi$ appartiennent à $H$.\\
\textit{2\up{d} critère d'équivalence relative}~: les diviseurs $\phi$ et $\psi$ sont équivalents si et seulement s'il existe deux diviseurs $\eta$, $\theta\in H$ tels que $\phi=\frac{\eta}{\theta}\psi$.\\
Kronecker esquisse une démonstration de la finitude du groupe des classes, en construisant un ensemble fini de diviseurs, tel que pour tout diviseur $\phi$, son inverse (dans le groupe des classes) soit l'un d'eux. \textit{Démonstration}. Soit $k$ un entier tel que~:
$$k^n\leq\text{Nm}\phi<(k+1)^n$$
On choisit une base de $O_K$, et on considère l'ensemble de tous les entiers algébriques à coefficients positifs $\leq k$. Il y en a $(k+1)^n>\text{Nm}\phi$, donc deux d'entre eux sont égaux modulo $\phi$, et leur différence est un diviseur principal $C\phi$, dont les coefficients sont inférieurs à $k$ en valeur absolue. Soit $C\phi=k\eta$ où $\eta$ est un nombre algébrique à coefficients inférieurs à 1. Alors $\text{Nm}\eta$ est borné par une grandeur $M$ indépendante de $\phi$, et 
$$\text{Nm}C.\text{Nm}\phi=k^n\text{Nm}\eta\leq\text{Nm}\phi.M$$
Donc $\text{Nm}C\leq M$, et l'inverse $C$ de $\phi$ peut être choisi dans l'ensemble fini des diviseurs de norme inférieure à $M$, \textit{q.~e.~d.}
\par
Kronecker avait vu dans les travaux de Gauss sur les formes quadratiques un cas particulier du <<~principe d'équivalence de Kummer~>>. Il existe en effet une relation étroite entre la théorie des entiers algébriques dans les extensions quadratiques de $\mathbb{Q}$ et la théorie des formes quadratiques de Gauss~: l'application $\alpha\mapsto\Fm\alpha$ associe à tout diviseur linéaire $x_1u_1+x_2u_2$ une forme quadratique en les indéterminées $u_1$ et $u_2$. Kronecker profite d'ailleurs pleinement de cette ressource historique pour prévenir cette <<~apparence d'étrangeté~>> que ne manquerait pas de ressentir son lecteur face à l'usage des indéterminées en théorie des nombres~:
\begin{quote}
Il suffit de rappeler l'introduction par Gauss des formes quadratiques en arithmétique pure, pour lever cette apparence d'étrangeté. Avant Gauss, on ne connaissait que des formes quadratiques de \textit{nombres}~; Gauss a le premier laissé tombé ce point de vue restreint suivant lequel on visait seulement la représentation des nombres par des formes quadratiques. Il a introduit, en arithmétique, des formes avec de véritables <<~indéterminées~>> (\textit{indeterminatae}).%
\footnote{\Cf \cite{kronecker1882} \S~22, p.~94-95.}
\end{quote}
\par
Les recherches décrites par Kronecker dans le reste de la section 19 constituent la préhistoire de la théorie du corps de classe. Helmut Hasse (\cf \cite{cassels1967}, ch.~XI, <<~History of Class Field Theory~>>) a comparé les contributions respectives de Kronecker, Weber et Hilbert à la genèse de cette théorie. Dedekind avait déjà remarqué%
\footnote{\Cf \cite{edwards1982}, p.~61, remarque n\up{o}33 de Dedekind.}, dans le cas d'un corps de nombres, qu'à cause de la finitude du groupe des classes, pour tout diviseur $\alpha$, il existe $m\in \mathbb{N}$ tel que $\alpha^m$ est principal, c'est-à-dire (à une unité près), égal à un élément $a\in O_K$~:
$$\alpha^m=\epsilon a$$
Dans l'anneau d'entiers de l'extension $K(\sqrt[m]{a})$, le diviseur $\alpha^m$ est donc égal, \textit{à une unité près}, à la puissance $m$-ième du diviseur principal $\sqrt[m]{a}$~:
$$\alpha^m=\epsilon (\sqrt[m]{a})^m$$
En décomposant en diviseurs premiers les deux membres de cette égalité, on en déduit que le diviseur $\alpha$ est égal à $\sqrt[m]{a}$ à une unité près, et est donc principal.
Comme il n'y a qu'un nombre fini de classes, il suffit d'adjoindre un nombre fini de radicaux au corps $K$ pour rendre principal chaque diviseur. Mais Kronecker ne se satisfait pas d'un tel procédé, qui n'offre pas une description suffisamment explicite des extensions de corps ainsi obtenues. Il explique qu'il avait, depuis l'hiver 1856, trouvé comment associer à certains corps de nombres une extension (un <<~genre associé~>>, \textit{associerte Gattung}) dont le degré est égal au nombre de classes, et telle que <<~toutes les propriétés profondes en rapport avec la composition et la répartition en classes du genre [...] se reflètent pour ainsi dire dans les propriétés élémentaires du genre associé~>>%
\footnote{\Cf \cite{kronecker1882}, p.~67.}. %
Pour mieux saisir la différence entre ces deux méthodes, soit $K=\mathbb{Q}(\sqrt{-23})$. On montre facilement que le groupe des classes est $\mathbb{Z}/3\mathbb{Z}$, et qu'il est engendré par la classe de $2u+\frac{1-\sqrt{-23}}{2}v$, qui est l'un des deux facteurs premiers de 2. La théorie hilbertienne du corps de classes consiste à déterminer une extension abélienne maximale non-ramifiée%
\footnote{\Cf \cite{cassels1967}, exercice 3. Plus précisément, l'extension abélienne en question doit être non-ramifiée et  <<~totalement décomposée en les valuations archimédiennes de $K$~>>. Dans le cas particulier ci-dessus, cela signifie qu'il y a exactement trois $K$-plongements de $L$ dans $\mathbb{C}$ essentiellement distincts vis-à-vis de la topologie induite par la norme dans $\mathbb{C}$.} %
$L$ de $K$, qui peut d'ailleurs s'obtenir, dans ce cas particulier, au moyen de la théorie des fonctions elliptiques. Il s'agit du corps de décomposition $L$ de l'équation $X^3-X-1=0$ sur $K$. Dans $L$, tous les diviseurs de $K$ deviennent principaux (mais $L$ peut très bien avoir de nouveaux diviseurs non principaux). Et $[L:K]$ est bien égal au nombre de classes de $K$, comme Kronecker l'affirme ($[L:K]=3$). Cette extension, non-ramifiée, fournit d'ailleurs l'exemple d'un <<~genre qui n'a pas de discriminant~>>, situation qui semblait susciter l'intérêt de Kronecker%
\footnote{\cf \cite{kronecker1882}, p.~22.}.  %
Démontrons-le directement. Soient $R$, $T$ et $S$ les anneaux d'entiers des corps $\mathbb{Q}$, $K$ et $L$. On a alors (\cite{cassels1967}, p.~17, \textit{tower formula})~:
$$\delta(S/R)=\delta(T/R)^{[L:K]}N_{K/\mathbb{Q}}\delta(S/T)$$
Or $T=\mathbb{Z}\left[1,\frac{1+i\sqrt{23}}{2}\right]$ et $\delta(T/R)=-23$. L'algorithme et la démonstration de Kronecker (\S~6 des \textit{Grundzüge}) pour construire une base du $Z$-module $O_L$, ont pour conséquence immédiate que $\delta(S/R)$ est inférieur en valeur absolue au discriminant des grandeurs $1, x_2, x_1, x_1x_2, x_1^2, x_1^2x_2$, où $x_1$ et $x_2$ désignent deux des racines de $X^3-X-1$. Or ce discriminant de 6 grandeurs n'est autre que $\mathcal{D}^3$, où $\mathcal{D}$ est le discriminant de $X^3-X-1$, comme Kronecker l'affirme au $\S~12$, et comme nous l'avons démontré plus haut. Le calcul donne donc~:
$$\vert\mathcal{D}^3\vert=23^3\geq\vert\delta(S/R)\vert=23^3\vert N_{K/\mathbb{Q}}\delta(S/T)\vert$$
donc $\vert\delta(S/T)\vert=1$, comme annoncé. Kronecker faisait certainement allusion à de tels exemples, tirés de la théorie des fonctions elliptiques. Sur le même exemple, l'autre méthode, critiquée par Kronecker, nous aurait conduit à calculer le diviseur
$$\left(2u+\frac{1-\sqrt{-23}}{2}v\right)^3=\epsilon\frac{3+i\sqrt{23}}{2}$$
où $\epsilon$ est une unité, puis à adjoindre au corps $K$ le radical $\sqrt[3]{\frac{3+i\sqrt{23}}{2}}$, rendant ainsi principal \mbox{$2u+\frac{1-\sqrt{-23}}{2}$}. Mais en procédant ainsi, on risque d'obtenir une extension ramifiée, et l'on est pas sûr d'obtenir une extension unique, indépendante du choix des radicaux adjoints au corps de base. L'histoire de la théorie du corps de classes a donc justifié la préférence de Kronecker pour l'autre méthode, bien que celui-ci n'ait pas pu, semble-t-il%
\footnote{\Cf les remarques de H. Hasse dans \cite{cassels1967}, ainsi que les remarques 32 à 36 de Dedekind et le commentaire de Edwards, Neumann et Purkert dans \cite{edwards1982}.}, %
s'élever des exemples tirés de la théorie des fonctions elliptiques à une théorie plus générale.

\section{Systèmes de diviseurs et codimension}
\label{Stufe}
Kronecker est donc parvenu à faire valoir de nouveau la théorie classique de la divisibilité (plus grand commun diviseur, décomposition en facteurs premiers) en étendant l'anneau d'entiers d'une extension algébrique par l'adjonction d'indéterminées. Au \S~20, il remarque que cette <<~conservation des déterminations conceptuelles lors du passage du rationnel à l'algébrique~>> ne vaut pas, par contre, lors du passage d'un corps ne contenant pas de variable à un corps contenant des variables (c'est-à-dire une extension transcendante d'un corps de nombres). Bien que la théorie développée dans les sections \S~14 à 17 s'applique aussi à de tels corps, Kronecker la juge à présent insuffisante. Le concept de diviseur devant exprimer <<~quelque chose de commun~>> (\textit{Gemeinsames}) à deux ou plusieurs grandeurs du domaine, il avait en effet été défini sur le modèle du <<~plus grand commun diviseur~>>, comme combinaison linéaire $x_1u_1+...+x_nu_n$ des grandeurs en question. Mais dans le cas où certaines des indéterminées de l'anneau d'entiers considéré sont conçues comme des \textit{variables}, le plus grand commun diviseur ne suffit plus à exprimer tout le <<~commun~>> de $n$ grandeurs. En effet, l'intuition géométrique fait correspondre au plus grand commun diviseur de $n$ polynômes une hypersurface, sur laquelle les $n$ polynômes s'annulent ensemble. Mais, en général, le lieu commun des zéros des $n$ polynômes peut contenir d'autres variétés, de dimensions inférieures. Le plus grand commun diviseur ne saurait en rendre compte. C'est donc l'intuition géométrique qui présente d'abord à notre conscience un <<~commun~>> autre que le <<~plus grand commun diviseur~>>, sous la forme d'une intersection, lieu des zéros communs à un ensemble de polynômes%
\footnote{L'opération intellectuelle par laquelle ce <<~quelque chose de commun~>>, concept d'abord vide de tout substrat mathématique, pourra légitimement subsumer le concept de plus grand commun diviseur et le concept d'intersection, n'est-il pas un bel exemple de synthèse (au sens kantien), renversant l'image d'un Kronecker conservateur obsédé par une vision arithmétisante de la réalité mathématique?}. %
Le concept de diviseur étant donc insuffisant pour rendre compte de ce <<~commun~>>, Kronecker définit à présent (\S~20--22) le concept de <<~système modulaire~>> (\textit{Modulsystem})%
\footnote{Ce terme est devenu \textit{modular system} chez Macaulay \cite{macaulay1916}.} %
ou <<~système de diviseurs~>> \label{KroneckerCycles} (\textit{Divisoren-System}). Un élément $G$ est dit <<~contenir le système de diviseurs $(F_1,...,F_n)$~>> lorsque $G \equiv 0 \text{~modulo~} (F_1,...,F_n)$. Mais comme pour échapper de nouveau à la notion ensembliste d'<<~idéal~>>, Kronecker représente tout système de diviseurs par une forme à coefficients indéterminés (par exemple $F_1u_1+...+F_nu_n$). Contrairement au cas des <<~diviseurs~>>, la relation <<~être contenu dans~>>, pour les systèmes de diviseurs, \textit{ne se traduit pas} par une relation de divisibilité.
\par
Kronecker introduit au \S~21 un concept de codimension (\textit{Stufe}) sans le définir de manière rigoureuse, en ayant recours encore une fois à l'intuition géométrique. On a vu que la méthode d'élimination du \S~10 permet de décomposer une variété affine, lieu des zéros d'un idéal de $k[x_1,...,x_n]$ sur $\overline{k}$, en composantes irréductibles. On peut alors définir la dimension comme dimension maximale de ses composantes. Au \S~20, il traite en détail le cas d'une variété affine de dimension 0 dans le cas de l'intersection complète, c'est-à-dire d'un ensemble fini de points solutions d'un système de $n$ équations à $(n-1)$ inconnues, sous l'hypothèse que le discriminant est non nul~; il utilise à cet effet la théorie du résultant de $n$ polynômes à $(n-1)$ indéterminées, sans d'ailleurs se soucier de savoir si la <<~résolvente~>> définie au \S~10 coïncide ou non avec le <<~résultant~>> classique%
\footnote{\label{ResolventeEtResultant}Macaulay a donné un contre-exemple, \cf \cite{macaulay1916}, p.~21.}. %
Dans ce cas particulier, la relation <<~être contenu dans~>> se traduit bien par une relation de divisibilité. Quant à l'origine de ces recherches, Kronecker renvoie à ses travaux de 1865 sur la formule d'interpolation de Lagrange. Nous avons déjà rencontré, au \S~4, l'usage qu'il fait de cette formule pour factoriser les polynômes. Dans un article de 1865, \textit{\ul{U}ber einige Interpolationsformeln für ganze Functionen mehrer Variabeln}, il généralise cette formule aux polynômes à plusieurs indéterminées. Le jacobien joue alors le rôle joué par la dérivée dans la formule de Lagrange. Soit un ensemble de $m$ points $\xi_k=(\xi_{1k},...,\xi_{nk}), 1\leq k\leq m$,  solution d'un système d'équations non homogènes à $n$ inconnues $F_1(x_1,...,x_n)=0, ..., F_n(x_1,...,x_n)=0$. On peut alors écrire (on le justifierait aujourd'hui au moyen du théorème des zéros de Hilbert)~:
$$F_i=(x_1-\xi_{1k})F_{1i}^{(k)}+...+(x_n-\xi_{nk})F_{ni}^{(k)}$$
Les déterminants $D_k(x_1,...,x_n)=\left\vert F_{hi}^{(k)}\right\vert_{1\leq h,i\leq n}$ permettent alors d'écrire la formule d'interpolation suivante qui définit une fonction $\mathfrak{F}(x_1,...,x_n)$ prenant la valeur $\mathfrak{F}_k$ au point $\xi_k$~:
$$\mathfrak{F}(x_1,...,x_n)=\sum_k\mathfrak{F}_k\frac{D_k(x_1,...,x_n)}{D_k(\xi_{1k},...,\xi_{nk})}$$
De plus, en remarquant qu'au point $\xi_k$
$$F_{hi}^{(k)}=\frac{\partial F_i}{\partial x_h},$$ 
on a~:
$$D_k(\xi_{1k},...,\xi_{nk})=J(\xi_{1k},...,\xi_{nk})$$
où $J$ est le déterminant fonctionnel $\left\vert\frac{\partial F_i}{\partial x_h}\right\vert_{1\leq h,i\leq n}$. En termes modernes, la formule d'interpolation de Kronecker décrit l'anneau de fonctions d'une variété affine intersection complète de dimension zéro, ensemble de $m$ points, comme étant de la forme $k^m$ (où $k$ est le corps des constantes), en associant à chaque fonction $\mathfrak{F}$ ses valeurs en ces points $(\mathfrak{F}_1,...,\mathfrak{F}_m)\in k^m$. 
\par
La forme dominante de la formule d'interpolation de Kronecker est égale à~:
$$\frac{R(x_1,..,x_n)}{\prod\deg F_i}\sum\frac{\mathfrak{F}(\xi_{1k},...,\xi_{nk})}{J(\xi_{1k},...,\xi_{nk})}$$
où $R$ est le jacobien des formes dominantes $f_1,...,f_n$ de $F_1,...,F_n$, soit $\left\vert\frac{\partial f_i}{\partial x_h}\right\vert_{1\leq h,i\leq n}$. Le facteur constant $\sum\frac{\mathfrak{F}(\xi_{1k},...,\xi_{nk})}{J(\xi_{1k},...,\xi_{nk})}$ intervient aussi dans un théorème de Jacobi (1836) que Kronecker démontre dans cette article~: si $\mathfrak{F}$ est de degré $< \deg J$, on a
$$\sum\frac{\mathfrak{F}(\xi_{1k},...,\xi_{nk})}{J(\xi_{1k},...,\xi_{nk})}=0.$$
La démonstration semble néanmoins incomplète. Kronecker affirme que la forme dominante de sa formule d'interpolation est, à un facteur constant près,
$$R(x_1,..,x_n)\sum\frac{\mathfrak{F}(\xi_{1k},...,\xi_{nk})}{J(\xi_{1k},...,\xi_{nk})}.$$
Il en conclut que, pour un tel $\mathfrak{F}$,
$$R(x_1,..,x_n)\sum\frac{\mathfrak{F}(\xi_{1k},...,\xi_{nk})}{J(\xi_{1k},...,\xi_{nk})}\in(f_1,...,f_n),$$
puis finalement que
$$R(x_1,..,x_n)\sum\frac{\mathfrak{F}(\xi_{1k},...,\xi_{nk})}{J(\xi_{1k},...,\xi_{nk})}=0.$$
Il ne justifie pas les deux dernières étapes de son raisonnement, bien qu'elles soient valables en général. En utilisant le formalisme de J.-P. Jouanolou, il suffit en effet de remarquer que l'idéal $(F_1,...,F_n)$ n'a pas de forme d'inertie de degré $\deg J$, et que l'idéal $(f_1,...,f_n)$ a exactement une forme d'inertie de degré $\deg J$ non triviale, le jacobien $R$.
\par
Le facteur constant $\sum\frac{\mathfrak{F}(\xi_{1k},...,\xi_{nk})}{J(\xi_{1k},...,\xi_{nk})}$ n'est autre que le <<~résidu de Grothendieck~>>%
\footnote{Serre le nomme <<~résidu de Grothendieck~>> dans \cite{serre1959},~p.81.}, %
que l'on exprime aujourd'hui comme trace d'un endomorphisme de multiplication (quand on travaille sur un corps algébriquement clos, la trace est bien la somme des valeurs d'une fonction en les différents points de la variété discrète). Pour mieux comprendre la situation, observons le cas d'une seule indéterminée. Le théorème de Jacobi dans ce cas était déjà connu d'Euler%
\footnote{\Cf tome 2 des \textit{Institutiones Calculi Integralis}, 1768--1770.}. %
Si $k$ est algébriquement clos, dans le $k$-module $k[X]/f(X)$, la trace d'une fonction $g$ est la somme des valeurs de $g$ en les racines de $f$. Les formules d'Euler disent que ($m=\deg f$)~:
$$\text{Tr}\left(\frac{x^i}{f'(x)}\right)=\begin{cases}0 & \text{si }0\leq i\leq m-2\\1 & \text{si }i=m-1\end{cases}$$
Et dans ce cas, la trace peut s'écrire au moyen d'un résidu de Cauchy~:
$$\text{\Huge\textgreek{e}}\left(\frac{x^i}{f(x)}\right)=\text{Tr}\left(\frac{x^i}{f'(x)}\right)$$
\par
Mais suivons à nouveau le fil des \textit{Grundzüge} pour observer comment Kronecker travaillait avec des variétés de dimension supérieure. Au \S~22, il expose une théorie des formes représentant les systèmes de diviseurs. Soit $A=k[x_1,...,x_{n-1}]$. Soit $S\subset A[u_1,u_2,...]$ l'ensemble (multiplicatif) des formes <<~proprement primitives~>>, c'est-à-dire dont les coefficients engendrent l'idéal $(1)$ de l'anneau $A$.
\begin{lemme}
Si les formes $E$ et $F$ sont proprement primitives alors la forme $EF$ l'est aussi, et réciproquement.
\end{lemme}
\textit{Démonstration}. Supposons $EF$ proprement primitive. Son idéal n'a donc pas de zéro dans $\overline{k}$. Or le produit des idéaux engendrés par les coefficients de $E$ et de $F$ contient l'idéal engendré par les coefficients de $EF$. Donc ce produit d'idéaux n'a pas de zéro non plus, donc $E$ et $F$ sont proprement primitives (théorème des zéros de Hilbert). Réciproquement, supposons $E$ et $F$ proprement primitives. Soit $z$ un zéro de l'idéal engendré par les coefficients de $EF$. Alors $EF(z)\equiv 0$, donc soit $E(z)\equiv 0$, soit $F(z)\equiv 0$. Absurde. Donc $EF$ est proprement primitive. \textit{q.~e.~d.}
\par
Kronecker travaille dans l'anneau $S^{-1}A[u_1,u_2,...]$. Une forme (en $u_1,u_2,...$) est dite de codimension $m$ quand la variété des zéros (dans $\overline{k}$) de ses coefficients (éléments de $A$) est de codimension $m$. Alors $m\leq n$, et $m=n$ si et seulement si la variété des zéros est vide, c'est-à-dire si la forme appartient à $S$ (Kronecker a ici recours à un cas du théorème des zéros de Hilbert). La forme $F_0$ <<~contient~>> la forme $F$ si et seulement si l'idéal des coefficients de $F_0$ est inclus dans l'idéal des coefficients de $F$. Mais cela \textit{ne se traduit pas} par une relation de divisibilité. Et cette relation <<~être contenu dans~>> donne lieu à une relation d'équivalence entre formes, qui ne se traduit pas non plus par l'égalité des formes dans $S^{-1}A[u_1,u_2,...]$. Le <<~plus grand commun contenu~>> (\textit{grössten gemeinsamen Inhalt}) de $F$ et $F_0$ est la forme $uF+F_0$ (où $u$ est une indéterminée n'entrant pas dans l'expression de $F$ et $F_0$). En termes de variétés, il s'agit de l'intersection des deux variétés.
\par
De même que l'adjonction d'indéterminées avait permis à Kronecker, aux \S~14--18, de transformer un anneau d'entiers en un anneau principal, il tente de l'utiliser, aux \S~20--22, pour transformer $A$ en un anneau $S^{-1}A[u_1,u_2,...]$, faisant de toute variété, de dimension donnée, une intersection complète (c'est-à-dire que l'idéal d'une variété de dimension $m$ devrait être engendré par $m$ polynômes). Afin de donner une formulation algébrique de la relation entre formes <<~être contenu dans~>> (qui n'est pas une relation de divisibilité), Kronecker énonce en effet le résultat suivant, pour deux <<~formes pures~>> de codimension $m$ (\textit{reine Formen m\up{ter} Stufe}) $F_0$, $F_1$~:
\begin{quote}
Soient $F_1,F_2,...,F_m$ des formes algébriques entières, qui ont chacune les mêmes coefficients et diffèrent donc entre elles seulement quant au système des indéterminées, alors l'une de ces formes est contenue dans une forme $F_0$ si $F_0$ est absolument équivalente [c'est-à-dire égale, à un élément de $S$ près] à une fonction linéaire homogène entière des $m$ formes $F_1,F_2,...,F_m$.%
\footnote{\cf \cite{kronecker1882}, \S~22, p.~91.}
\end{quote}
Kronecker ajoute~:
\begin{quote}
  Ceci peut servir de définition de la codimension (\textit{Stufenzahl}) $m$ d'une forme pure $F$, si l'on ajoute seulement qu'aucun nombre plus petit que $m$ ne peut servir à une telle représentation.%
\footnote{\Cf \cite{kronecker1882}, p.~92.}%
\end{quote}
Il nous est difficile de juger ces deux affirmations, qui lui ont probablement été suggérées par l'étude des variétés affines \textit{linéaires} (dans ce cas l'énoncé ci-dessus est évidemment juste). Dans le cas général, les travaux récents de Scheja-Storch \cite{scheja2001} ne laissent pas entrevoir de relation évidente entre ce plus petit entier $m$ et la codimension (au sens géométrique du terme). Une définition algébrique de la codimension est bien ce qui manquait à Kronecker, pour pouvoir étendre sa théorie des <<~systèmes de diviseurs~>> à des anneaux $A$ plus généraux. La théorie de l'élimination sur laquelle il fonde, au \S~10, la notion de \textit{Stufe} se rapporte exclusivement à des anneaux de polynômes \textit{sur un corps}, et à l'intuition géométrique sous-jacente des variétés affines. Il faudrait faire une théorie de l'élimination sur $\mathbb{Z}$ pour dépasser ce cadre. Kronecker semblait convaincu qu'il était possible de franchir le pas. Au \S~21, il décrit à quoi ressemblerait la théorie de la dimension dans le cas $A=\mathbb{Z}[x_1,x_2,x_3]$~:
\begin{quote}
  Prenons par exemple $n=4$ et pensons les trois variables comme des coordonnées quelconques de l'espace, alors les diviseurs de codimension 1 sont soit des nombres soit des fonctions entières des coordonnées, dont l'annulation représente des surfaces. Parmi les systèmes de diviseurs de codimension 2, il s'en trouve dont les éléments ne peuvent tous s'annuler simultanément, et dont l'un des éléments peut être choisi comme nombre, mais il s'en trouve aussi, dont l'annulation simultanée des éléments représente une courbe~; parmi les systèmes modulaires de codimension 3, l'on en trouve qui représente de même des systèmes de points...%
\footnote{\Cf \cite{kronecker1882}, \S~21.}%
\end{quote}
Il indique même que la théorie des systèmes de diviseurs de codimension 2 dans $\mathbb{Z}[x]$ produit une théorie des nombres algébriques, la décomposition du système de diviseurs $(F(x),p)$ correspondant à la décomposition de $p$ en facteurs premiers dans l'anneau d'entiers algébriques sur $\mathbb{Q}$ de $\mathbb{Q}[X]/F(X)$ (sauf pour un ensemble fini de valeurs de $p$). En fait, il faut attendre la théorie des schémas de Grothendieck pour avoir, comme le souhaitait Kronecker, une théorie des variétés (les <<~schémas~>>) sur $\mathbb{Z}$ ou un anneau quelconque.
\par
Quant à la <<~grandeur entière~>>, il y avait en tout cas une lacune dans les \textit{Grundzüge}~: bien que le concept de grandeur entière et celui de plus grand commun diviseur aient subi avec succès le passage à des extensions algébriques quelconques de $\mathbb{Q}$, voire de $\mathbb{Q}(T)$, dans les sections \S~14 à 18, et bien qu'au concept de plus grand commun diviseur ait été substitué celui d'intersection dans les sections \S~20 à 22, effectuant ainsi le passage d'un corps de nombres $k$ à des extensions purement transcendantes de dimension quelconque $k(x_1,...,x_{n-1})$ au moyen d'une théorie de l'élimination dans $k[x_1,...,x_{n-1}]$ qui rend compte des variétés de codimension supérieure à 1, il restait à faire une théorie de la <<~grandeur entière~>> dans les extensions algébriques de $k(x_1,...,x_{n-1})$, qui tienne compte elle aussi des variétés de codimension supérieure à 1. Qu'il nous soit permis d'employer un vocabulaire moderne~: il fallait en quelque sorte définir et étudier des <<~morphismes finis~>> entre variétés affines.
\par
Un court article \cite{kronecker1883} publié en 1883 témoigne encore du caractère non achevé des \textit{Grundzüge} et montre que Kronecker était à la recherche d'une autre formulation algébrique de la relation entre systèmes de diviseurs <<~être contenu dans~>>. Il y démontre le théorème suivant~:%
\footnote{\cf \cite{kronecker1883}, p.~419--420. Nous avons corrigé une faute dans la dernière formule de cette citation (le texte des \textit{Werke} contient $r$ au lieu de $\rho-r$).}%
\begin{quote}
Soient
$$M_0, M_1, M_2,..., M_{n+1}$$
des grandeurs entières d'un domaine de rationalité $[\rat',\rat'',\rat''',...]$, alors l'équation~:
$$\sum_hM_hU^h.\sum_iM_{m+i}U^{i-1}=\sum_kM_k^\prime U^k
\quad\left(\begin{array}{l}
h=0,1,...m\\
i=1,2,...n-m+1\\
k=0,1,...n
\end{array}\right),$$
où $U$ est une indéterminée, définit $n+1$ grandeurs $M_0^\prime, M_1^\prime,... M_n^\prime$ du même domaine $[\rat',\rat'',\rat''',...]$.\\
(...)\\
Le produit de formes~:
$$\sum_hM_hU_h.\sum_iM_{m+i}U_{m+i}
\quad\left(\begin{array}{l}
h=0,1,...m\\
i=1,2,...n-m+1
\end{array}\right)$$
vérifie une équation algébrique de degré $\rho$, dans laquelle le coefficient de la puissance $\rho$\up{ième} est égal à 1, mais celui de la puissance $r$\up{ième} \textit{est une forme contenant le produit de formes}
$$\prod_h\sum_kM_k^\prime V_{hk}
\quad\left(\begin{array}{l}
h=1,2,...\rho-r\\
k=0,1,...n
\end{array}\right)$$
\end{quote}
La relation <<~être contenu dans~>> qu'utilise Kronecker dans l'énoncé de ce théorème est simplement la relation d'inclusion entre les idéaux engendrés par les coefficients des formes. Autrement dit, le coefficient de la puissance $r$\up{ième} ($0\leq r\leq\rho$) appartient à l'idéal $(M_0^\prime, M_1^\prime,..., M_n^\prime)^{\rho-r}$. A la suite de ce théorème, Kronecker définit une nouvelle relation <<~être contenu dans~>> (plus faible que l'inclusion entre idéaux, mais qui implique au moins l'inclusion entre les radicaux des idéaux). Il explique ensuite, dans un paragraphe quelque peu obscur, que cette nouvelle relation permet de définir un nouveau concept de <<~grandeur entière~>>%
\footnote{Comme nous l'avons mentionné, Kronecker souhaitait réduire la théorie des nombres algébriques à celle des systèmes de diviseurs de rang 2 dans $\mathbb{Z}[x]$, et plus généralement toute la théorie des <<~formes entières algébriques de rang $\mathfrak{m}$~>> à une théorie des <<~formes entières rationnelles de rang $\mathfrak{m}+1$~>> (\cf \cite{kronecker1882} p.~113).}. %
Cet article, son titre l'indique, se rapportait principalement aux variétés de codimension supérieure à 1. Sans recours à l'historiographie, il nous serait impossible de comprendre ce que Kronecker visait. 
\par
J. Molk, dans un long commentaire en français des \textit{Grundzüge}, nous offre une explication possible quand il se propose de décomposer les systèmes de diviseurs en produits de systèmes irréductibles%
\footnote{\Cf \cite{molk1885} \S~IV.1. Il faut ici entendre la notion de produit, avec Molk, au sens du produit des idéaux correspondant aux systèmes de diviseurs. Molk affirme (cf. \cite{molk1885} p.~107) que sa méthode, décrite pour $A=k[x,y]$, est aussi valable pour $A=\mathbb{Z}[x]$.}, %
pour $A=k[x,y]$. A cet effet, dans le cas d'un système de diviseurs $(\Phi(x,y),\Psi(x,y))$, il montre que ce système est équivalent au système des coefficients de la résolvente de $\Phi$ et $\Psi$ (\cf \cite{molk1885} p.~79-97). Il en est de même en général pour un système de diviseurs $(\Lambda_1(x,y),\Lambda_2(x,y),...,\Lambda_\mu(x,y))$, \textit{à condition} d'affaiblir la relation d'équivalence entre systèmes de diviseurs comme le faisait Kronecker en 1883 (\cf \cite{molk1885} p.~104-105). Malgré ses efforts, Molk semble abandonner tout espoir d'une théorie générale des systèmes de diviseurs qui rende compte des facteurs multiples dans la résolvente. Or la relation d'équivalence entre systèmes de diviseurs consistant en l'égalité des radicaux des idéaux donne une théorie générale sans notion de multiplicité~; mais Kronecker et Molk cherchaient une relation d'équivalence plus forte. La recherche d'une telle relation d'équivalence était probablement motivée par la recherche d'une notion de multiplicité, et renoncer à toute notion de multiplicité était donc un constat d'échec. 
\par
D'autres commentateurs des \textit{Grundzüge} (Hurwitz, H. M. Edwards) ont vu dans l'article de 1883 le chaînon manquant qui permet au moins de démontrer rigoureusement les considérations des sections 14 à 18. Il y avait en effet, comme nous l'avons mentionné, un point obscur dans la section 14, qui a fait dire à Dedekind~:
\begin{quote}
  Ici se cache le c{\oe}ur de la théorie des idéaux de Kronecker~; la simplicité qui nous étonnait au début apparaît bientôt sous une autre lumière, si l'on exige que la démonstration soit complètement achevée.%
\footnote{\Cf \cite{edwards1982}, remarque n\up{o}20 de Dedekind, et le commentaire des éditeurs que nous résumons ci-dessous.}.
\end{quote}
\label{ChainonManquant}
Kronecker semble en effet déduire, du fait qu'une forme $Q$ vérifie une équation $\text{Nm}(X-Q)=0$ à coefficients entiers sur le corps de base, le fait que cette forme soit elle-même à coefficients entiers.
Dans le théorème de l'article de 1883, si l'on spécialise les indéterminées $U_i$ ($0\leq i\leq n$), on obtient facilement le corollaire suivant~:
\begin{corol}
  Les grandeurs $M_kM_{m+i}$ (pour tout $k$, $i$ tels que $0\leq k\leq m$ et $1\leq i\leq n-m+1$) sont entières sur $\mathbb{Z}[M_0^\prime, M_1^\prime,..., M_n^\prime]$.
\end{corol}
Ce corollaire se généralise facilement au produit de formes à un nombre quelconque d'indéterminées~:
$$\sum_hM_h\phi_h.\sum_iM_{m+i}\psi_i=\sum_kM_k^\prime\chi_k$$ 
où $\phi_h$, $\psi_i$, $\chi_k$ désignent des produits quelconques des indéterminées%
\footnote{\Cf par exemple \cite{edwards1990}, partie 0, corollaire 3 p.~4. \note{ou Konig1904 p.~78? (référence donnée par macaulay1916 p.~107)}}. %
En prenant $\sum_hM_h\phi_h=X-Q$ et $\sum_iM_{m+i}\psi_i=\frac{\text{Nm}(X-Q)}{X-Q}$, on en déduit la proposition annoncée.
\\\par
Reprenons le fil des \textit{Grundzüge}. Kronecker entendait réduire la théorie des <<~formes entières algébriques de rang $\mathfrak{m}$~>> à une théorie des <<~formes entières rationnelles de rang $\mathfrak{m+1}$~>>. Dans la vingt-cinquième et dernière section des \textit{Grundzüge}, il tente d'appliquer ce principe pour $\mathfrak{m}=1$~; il donne ainsi un nouvel algorithme (p.~114--117), déjà annoncé ailleurs dans les \textit{Grundzüge}, pour la décomposition d'une grandeur algébrique entière en ses diviseurs irréductibles. Dans le cas $O_k=\mathbb{Z}$, l'idée de cette méthode n'est pas nouvelle\note{Schönemann, \cf p.~118 et 80, ou bien Dedekind, \cf p.~80}. Plus généralement, soit
$$\begin{array}{l}
O_k=\mathbb{Q}[\rat',\rat'',...]\\
k=\mathbb{Q}(\rat',\rat'',...)\\
K=\mathbb{Q}(\mathfrak{G},\rat',\rat'',...)
\end{array}$$
où $\mathfrak{G}$ vérifie une équation%
\footnote{Nous ne décrivons ici que le cas où $K$ est une extension algébrique de $k$ ayant un élément primitif $\mathfrak{G}$. Mais on peut toujours se ramener à cette situation en adjoignant des indéterminées $u', u'',...$ ; on aurait alors $K=\mathbb{Q}(u',u'',...)(\mathfrak{G},\rat',\rat'',...)$ et cela ne nuit pas à la suite du raisonnement. Selon Edwards, c'est Hensel\note{Crelle 113(1894):61--83} qui a démontré, dans le cas particulier $k=\mathbb{Q}$, la validité de l'algorithme décrit ici par Kronecker. Edwards avoue <<~\textit{I do not know of a proof of Kronecker's more general case. (If one could be given, it would be a large step toward the solution of the problem of factoring divisors in the general case.).}~>>. Nous essayons ici de reconstituer le raisonnement, imparfait, de Kronecker dans le cas où $k=\mathbb{Q}(\rat',\rat'',...)$ et $O_k=\mathbb{Q}[\rat',\rat'',...]$.}
$$\mathfrak{F}(\mathfrak{G})=0$$
avec $\mathfrak{F}(\rat)\in\mathbb{Q}[\rat',\rat'',...][\rat]$. 
\par
Ici Kronecker ne parvient pas à exprimer une base agréable de $O_K$ avec $\mathfrak{G}$ et ses puissances (p.~111-112). En tout cas, $O_K$ contient l'anneau engendré par $\mathfrak{G}, \rat', \rat'',...$~:
$$\mathbb{Q}[\mathfrak{G},\rat',\rat'',...]\subset O_K$$
Nous allons préciser cette affirmation. Remarquons que $K$, anneau engendré par $\mathfrak{G}$ sur $k$, est isomorphe à un anneau quotient~:
$$K\simeq k[\rat]/\mathfrak{F}(\rat)$$
Soit $\Delta$ le discriminant de $\mathfrak{F}(\rat)$. Alors $\Delta\in O_k$. On va montrer que
$$O_K\subset O_{k,(P)}[\mathfrak{G}]$$
où l'anneau de droite est engendré par $\mathfrak{G}$ sur $O_{k,(P)}$, localisé de $O_k$ par rapport à n'importe quel idéal premier $(P)$ avec $P$ ne divisant pas $\Delta$. Soit en effet $x\in O_K$, alors il existe des $a_i\in k$ tels que~:
$$x=\sum_{i=0}^{\deg\mathfrak{F}-1}a_i\mathfrak{G}^i$$
Les $(\deg\mathfrak{F})$ équations conjuguées forment un système d'équations linéaires en les $a_i$, dont le déterminant est un déterminant de Vandermonde ayant pour carré $\Delta$. Les formules de Cramer permettent d'exprimer chaque $a_i$ comme quotient d'un élément de $O_k[\mathfrak{G}]$ par ce déterminant. Si $P\in O_k$ ne divise pas $\Delta$, il est donc clair que chaque $a_i$ appartient à l'anneau $O_{k,(P)}[\mathfrak{G}]$. \textit{q.~e.~d.}
\par
Soit une telle grandeur entière $P$ irréductible dans $O_k$ et ne divisant pas $\Delta$. Essayons de la décomposer en produit de diviseurs (rappelons, conformément à la théorie des diviseurs exposée plus haut, qu'il s'agit de décomposer $P$ dans $O_K\otimes S^{-1}O_k[u_1,u_2,...]$). Remarquons tout d'abord que l'on a bien une théorie de la codimension pour les systèmes de diviseurs dans $\mathbb{Q}[\rat',\rat'',...][\rat]$, théorie garantie par la méthode d'élimination de Kronecker et son interprétation géométrique%
\footnote{On serait gêné par l'absence d'une théorie de la codimension dans le cas $O_k=\mathbb{Z}[\rat',\rat'',...]$.}. %
Or pour décomposer $(P)$ en ses composantes de codimension 1 (c'est-à-dire en éléments de $O_K$), on va décomposer $(P,\mathfrak{F}(\rat))$ en ses composantes de codimension 2 (c'est-à-dire que l'on cherche un produit d'idéaux qui, en codimension 2, présente le même lieu des zéros que $(P,\mathfrak{F}(\rat))$). A cet effet, on a naturellement recours à une décomposition de $\mathfrak{F}$~:
$$\mathfrak{F}(\rat)=f_1(\rat)^{n_1}f_2(\rat)^{n_2}...+P\overline{\mathfrak{F}}(\rat)$$
avec $f_1(\rat), f_2(\rat),...$ irréductibles modulo $P$ (cela a bien un sens puisque l'anneau $O_{k,(P)}/(P)$ est un corps et que l'anneau de polynômes $\left(O_{k,(P)}/(P)\right)[\rat]$ est donc factoriel). D'ailleurs, comme le discriminant $\Delta$ de $\mathfrak{F}$ est non nul modulo $P$, une telle décomposition ne compte pas de facteur double et $n_1=n_2=...=1$.
\par
\textit{Remarque~:} une telle décomposition de $\mathfrak{F}$ donne bien une décomposition <<~géométrique~>> de $(P,\mathfrak{F}(\rat))$, au sens suivant. Mettons toutes les fractions rationnelles au même dénominateur~:
$$f_1(\rat)=\frac{\tilde{f_1}(\rat)}{c(\rat',\rat'',...)},\quad f_2(\rat)=\frac{\tilde{f_2}(\rat)}{c(\rat',\rat'',...)},...\quad \overline{\mathfrak{F}}(\rat)=\frac{\tilde{\mathfrak{F}}(\rat)}{c(\rat',\rat'',...)}$$
avec $c(\rat',\rat'',...)$ non divisible par $P$, de sorte que
$$c(\rat',\rat'',...)\mathfrak{F}(\rat)=\tilde{f_1}(\rat)\tilde{f_2}(\rat)...+P\tilde{\mathfrak{F}}(\rat)$$
Les lieux géométriques suivants sont alors confondus~:
$$\begin{array}{l@{}c@{}l}
Z\left(\prod\limits_i\left(P,\tilde{f_i}(\rat)\right)\right)&=&Z\left(P,\prod\limits_i\tilde{f_i}(\rat)\right)\\
&=&Z(P,c(\rat',\rat'',...)\mathfrak{F}(\rat))\\
&=&Z\left((P,c(\rat',\rat'',...))(P,\mathfrak{F})\right)
\end{array}$$
D'où 
$$\bigcup\limits_iZ\left(P,\tilde{f_i}(\rat)\right)=Z\left(P,c(\rat',\rat'',...)\right)\cup Z\left(P,\mathfrak{F}\right)$$
\par
Revenons à la théorie des diviseurs. La décomposition de $\mathfrak{F}(\rat)$ donne, en posant $\rat=\mathfrak{G}$~:
$$\tilde{f_1}(\mathfrak{G})\tilde{f_2}(\mathfrak{G})...+P\tilde{\mathfrak{F}}(\mathfrak{G})=0$$
Donc dans $O_K$, $P$ divise $\prod\tilde{f_i}(\mathfrak{G})$. Donc dans $O_K\otimes S^{-1}O_k[u_1,u_2,...]$, il divise
$$\prod\left(P+u_i\tilde{f_i}(\mathfrak{G})\right)$$
Il est clair d'autre part que $P+u_i\tilde{f_i}(\mathfrak{G})$ divise $P$ pour tout $i$.
\begin{lemme}
 Pour tous $i\neq j$, le diviseur $P+u_i\tilde{f_i}(\mathfrak{G})+u_j\tilde{f_j}(\mathfrak{G})$ est une unité.
\end{lemme}
\textit{Démonstration~:} on écrit l'égalité de Bézout suivante dans l'anneau $\left(O_{k,(P)}/(P)\right)[\rat]$~:
$$A(\rat)f_i(\rat)+B(\rat)f_j(\rat)=C$$
Et $C\neq 0$ car $f_i$ et $f_j$ sont premiers entre eux. On relève cette égalité dans l'anneau $O_{k,(P)}[\rat]$~:
$$\tilde{A}(\rat)\tilde{f_i}(\rat)+\tilde{B}(\rat)\tilde{f_j}(\rat)=\tilde{C}$$
avec $\tilde{C}$ non divisible par $P$. On pose dans cette égalité $\rat=\mathfrak{G}$~; comme le diviseur $u_i\tilde{f_i}(\mathfrak{G})+u_j\tilde{f_j}(\mathfrak{G})$ divise $\tilde{f_i}(\mathfrak{G})$ et $\tilde{f_j}(\mathfrak{G})$, il divise $\tilde{C}$ aussi. Donc le diviseur $u_i\tilde{f_i}(\mathfrak{G})+u_j\tilde{f_j}(\mathfrak{G})+P$ divise le diviseur $P+u\tilde{C}$ qui lui-même est une unité, donc celui-là en est une aussi, \textit{q.~e.~d.}
\par
Ce lemme montre que les $P+u_i\tilde{f_i}(\mathfrak{G})$ sont tous premiers entre eux. Mais donc comme chacun divise $P$, leur produit aussi. Finalement les diviseurs $P$ et $\prod\left(P+u_i\tilde{f_i}(\mathfrak{G})\right)$ se divisent mutuellement et sont égaux, à une unité près. On a donc trouvé une décomposition de $P$. Il reste à vérifier que les diviseurs $P+u_i\tilde{f_i}(\mathfrak{G})$ sont irréductibles.
\par
\begin{lemme}
Pour tout $i$, le diviseur $P+u_i\tilde{f_i}(\mathfrak{G})$ est irréductible.
\end{lemme}
\textit{Démonstration~:} supposons qu'il existe $f(\mathfrak{G})\in O_K$ tel que $P+u_i\tilde{f_i}(\mathfrak{G})+uf(\mathfrak{G})$ soit différent d'une unité, et montrons qu'alors $P+u_i\tilde{f_i}(\mathfrak{G})$ divise $f(\mathfrak{G})$. Ecrivons une égalité de Bézout dans $\left(O_{k,(P)}/(P)\right)[\rat]$ que l'on relève immédiatement à $O_{k,(P)}[\rat]$~:
$$A(\rat)\tilde{f_i}(\rat)+B(\rat)f(\rat)=C+P.D(\rat)$$
Quitte à réduire au même dénominateur $A(\rat)$, $B(\rat)$, $C$ et $D(\rat)$, on peut même supposer qu'ils sont tous éléments de $O_k[\rat]$. Posons $\rat=\mathfrak{G}$, on obtient une égalité de la forme~:
$$A(\mathfrak{G})\tilde{f_i}(\mathfrak{G})+B(\mathfrak{G})f(\mathfrak{G})=C+P.D(\mathfrak{G})$$
avec $A(\mathfrak{G})$, $B(\mathfrak{G})$, $C$ et $D(\mathfrak{G})$ éléments de $O_K$. On en déduit que le diviseur $u_i\tilde{f_i}(\mathfrak{G})+uf(\mathfrak{G})+P$ divise le diviseur $C+vP$. Ce dernier ne peut donc pas être une unité, et if faut que $P$ divise $C$. Alors $f_i(\rat)$ et $f(\rat)$ ont un facteur commun modulo $P$ (c'est-à-dire si l'on considère ces deux polynômes comme éléments de l'anneau $\left(O_{k,(P)}/(P)\right)[\rat]$). Mais $f_i$ est irréductible modulo $P$, donc il divise $f$ modulo $P$. Cette relation de divisibilité, relevée à $O_{k,(P)}[\rat]$ et quitte à réduire toutes les fractions au même dénominateur, implique que dans $O_k[\rat]$, il existe un $C\in O_k$ non divisible par $P$ tel que
$$C.f(\rat)\in(\tilde{f_i}(\rat),P)$$
D'où, en posant $\rat=\mathfrak{G}$,
$$P+u_i\tilde{f_i}(\mathfrak{G})\quad\vert\quad C.f(\mathfrak{G})$$
et comme $C$ n'est pas divisible par $P$, on a finalement
$$P+u_i\tilde{f_i}(\mathfrak{G})\quad\vert\quad f(\mathfrak{G}),$$
ce qu'il fallait démontrer.
\\\par
Retenons seulement que le texte des \textit{Grundzüge} est de nature programmatique. Le titre l'annonçait déjà, <<~Esquisse d'une théorie arithmétique des grandeurs algébriques~>>. Kronecker n'hésite pas à sauter les points obscurs et laisser les détails à l'attention des commentateurs et de la recherche future, dans un ouvrage qui est le fruit d'une longue recherche menée sur des cas particuliers, champs d'essai des algorithmes, un ouvrage qui n'est pas une démonstration \textit{a priori} de résultats pressentis. Il n'a pourtant pas négligé de chercher les fondements simples de sa théorie (\textit{einfachste Grundlagen}, comme il le dit dans l'article de 1883)~; il définit certes des concepts (ceux de <<~genre~>>, de <<~grandeur entière~>>, de discriminant, de <<~système de diviseurs~>>, de codimension,...). Mais il ne semble pas vouloir en donner une construction logique suffisante. L'ouvrage s'organise alors autour d'algorithmes qui mettent en {\oe}uvre ces concepts et en donnent une certaine légitimation (l'algorithme de factorisation au \S~4, l'algorithme de calcul d'une base du $O_k$-module $O_K$ au \S~6, l'algorithme pour le calcul de la \textit{Gesammtresolvente} au \S~10, l'algorithme de décomposition en diviseurs premiers aux \S~18 et 25...), et autour de cas particuliers (celui du corps de décomposition d'une équation aux \S~11 et 12, celui des corps de nombres aux \S~14 à 19, celui des variétés de dimension 0 au \S~20...) et de l'évocation de travaux antérieurs (par exemple ses recherches sur le corps de classes et la multiplication complexe au \S~19, ou celles sur les extensions non ramifiées -- les <<~genres sans discriminant~>>~-- au \S~8), cas particuliers qui entraînent peu à peu l'intuition vers son nouveau vêtement.

\end{document}